\documentclass[12pt,a4paper,intlimits,sumlimits]{amsart}

\usepackage[usenames,dvipsnames,svgnames,table,rgb]{xcolor}
\usepackage{graphicx,txfonts,latexsym,fancyhdr,calc}
\usepackage{amsfonts,amsmath,amsthm,amssymb}
\usepackage[
allcolors=blue,allbordercolors=blue,pdfborderstyle={/S/U/W 1}]{hyperref}
\usepackage[ocgcolorlinks]{ocgx2}
\usepackage{epigraph}
\usepackage{tcolorbox}
\usepackage{asymptote}
\usepackage{parskip}
\usepackage{pgfgantt}
\usepackage{tcolorbox}
\usepackage{xfrac}

\ifdefined\SMART 
\usepackage[paperwidth=12cm,paperheight=24cm,left=5mm,right=5mm,top=10mm,bottom=5mm]{geometry}
\else
\usepackage[left=2cm,right=2cm,bottom=2cm,top=2cm]{geometry} 
\fi

\usepackage[normalem]{ulem}

\usepackage{txfonts,pxfonts,tikz} 
\usepackage{tgschola}
\usepackage[T1]{fontenc}

\theoremstyle{plain}
\newtheorem{theorem}{Theorem}

\newtheorem{lemma}{Lemma}
\newtheorem{corollary}[theorem]{Corollary}
\newtheorem{proposition}{Proposition}

 \newtheorem*{reptheorem}{Theorem \ref{th:op}}

\theoremstyle{definition}

\newtheorem{remark}{Remark}

\newtheorem{que}{Question}

\newtheorem{case[theorem]}{Case}

\newcommand{\beql}[1]{\begin{equation}\label{#1}}
\newcommand{\eeq}{\end{equation}}

\newcommand{\Abs}[1]{{\left|{#1}\right|}}

\newcommand{\Set}[1]{{\left\{{#1}\right\}}}

\newcommand{\RR}{{\mathbb R}}

\newcommand{\ZZ}{{\mathbb Z}}

\newcommand{\supp}{{\rm supp\,}}

\newcommand{\ft}[1]{\widehat{#1}}

\newcounter{rem}
\newcounter{step}
\newcounter{mysec}
\newcounter{mysubsec}[mysec]
\title{Spectra for finite unions of line segments}

\author{Mihail N. Kolountzakis}
\address{\href{https://math.uoc.gr/en}{Department of Mathematics and Applied Mathematics}, University of Crete,\\Voutes Campus, 70013 Heraklion, Greece,\newline and \newline \href{https://ics.forth.gr/}{Institute of Computer Science}, Foundation of Research and Technology Hellas, N. Plastira 100, Vassilika Vouton, 700 13, Heraklion, Greece}
\email{kolount@gmail.com}
\thanks{M. Kolountzakis acknowledges the generous hospitality and support of Fudan University.}

\author{Ruxi Shi}
\address{\href{https://scms.fudan.edu.cn/}{Shanghai Center for Mathematical Sciences, Fudan University}, 200438
Shanghai, China}
\email{ruxishi@fudan.edu.cn}
\thanks{R. Shi is supported by NSFC No. 12571198, NSFC No. 12231013 and the New Cornerstone Science Foundation through the New Cornerstone Investigator Program.}

\author{Sha Wu}
\address{\href{https://math.uoc.gr/en}{Department of Mathematics and Applied Mathematics}, University of Crete,\\Voutes Campus, 70013 Heraklion, Greece.\newline and \newline
\href{https://math.hnu.edu.cn/index.htm}{School of Mathematics, Hunan University}, 410082 Changsha, China}
\email{shaw0821@163.com}
\thanks{S. Wu is supported by Hunan Provincial Innovation Foundation for Postgraduate (LXBZZ2024024) and  grateful for the financial support provided by the China Scholarship Council (CSC)}

\begin{document}
\makeatletter
\@namedef{subjclassname@2020}{\textup{2020} Mathematics Subject Classification}
\makeatother
\date{22 November 2025}

\begin{abstract}
In this paper we study the spectrality of arc-length measures supported on the union of two line segments in the plane. We show that any such spectral measure must admit a line spectrum. Moreover, when the two segments are non-parallel, such   spectral measure admits only line spectra. Thus, in this case every spectrum is one dimensional.  In addition we show that this property fails for unions of three or more segments in the plane. We construct some arc-length spectral measures supported on the union of at least three line segments such that none of its spectra is contained in a line. Finally, we work in the general framework of arc-length  measures supported on 
finite unions of curves in $\mathbb{R}^d$. We show that the size of any orthogonal set for such a measure inside a ball of radius $R$ grows at most linearly in $R$.   We also give an alternative proof of this bound, and in fact obtain a more general result of growth  rate of orthogonal sets for Ahlfors--David regular measures in $\mathbb{R}^d$ (not restricted to the one-dimensional setting).
\end{abstract}
 \maketitle
\tableofcontents

\section{Introduction}

\subsection{Line  spectra for arc-length measures supported on two line segments}
A finite Borel measure 
$\mu$ on $\RR^{d}$
 is called a {\it spectral measure} if there exists a countable set $\Lambda\subset\RR^d$
 such that the exponential functions
$E_\Lambda=\{e^{2\pi i<\lambda, x>}: \lambda\in\Lambda\}$
form an orthonormal basis of $L^{2}(\mu)$.
Such a set  $ \Lambda$ is called a {\it spectrum} of $\mu$. In particular, if the  Lebesgue measure restricted on the set $\Omega$ is  a spectral measure,  then we say $\Omega$ is a {\it spectral set}.  The theory of spectral sets traces back to the work of Fuglede \cite{fuglede1974operators} who proposed 

{\bf{Fuglede Conjecture:}} $\Omega\subset \mathbb{R}^d$ is a  spectral set  if and only if  it tiles $\mathbb{R}^d$  by translations.

Both directions of the conjecture fail when $d\geq 3$  \cite{tao2004fuglede,kolountzakis2006tiles,kolountzakis2006hadamard,farkas2006onfuglede,farkas2006tiles}, but  the conjecture is still open in both directions when $d=1,2$. Moreover, the problem of deciding whether a given measure  is spectral, and of describing the possible spectra, is closely related to the Fuglede conjecture and has been extensively studied in recent decades for Lebesgue measure, for self-similar measures, and for various classes of singular measures. 

Spectra can be surprisingly rich. A particularly natural class of spectra to consider  for measures in $\RR^d$  are line spectra, that is, spectra contained in a straight line. One may then ask to what extent the geometry of the support forces spectra to be one dimensional. In this paper we focus on the case where the measure is supported on the union of finitely many line segments in the plane   and is equal to a constant multiple of arc-length on each of these segments. We are interested in whether the spectra of such arc-length spectral measures lie on straight lines.
In particular, when the arc-length measure is supported on the union of two line segments, all known results so far show that every spectrum of the associated spectral measure is contained in   straight lines. For example, several recent works have investigated the spectrality of the   measure
\[
\mu=\frac{1}{2}L_{[t, t+1]}\times\delta_0+\frac{1}{2}\delta_0\times L_{[t, t+1]},
\]
 where $L_{[t, t+1]}$ is the Lebesgue measure restricted on the interval $[t, t+1]$ and $\delta_0$ is the Dirac measure on the origin. 
Combining these papers \cite{kolountzakis2025non,Lev18,lai2021spectral,ai2023spectrality} with the  paper \cite{kolountzakis2025spectrality} by the first and third author, one obtains necessary and sufficient conditions for the measure $\mu$  to be a spectral measure and shows that it has only line spectra. See also the relevant works \cite{kolountzakis2025non,chen2025fourier}

\begin{figure}[htbp]
  \centering
\includegraphics[width=1.0\textwidth]{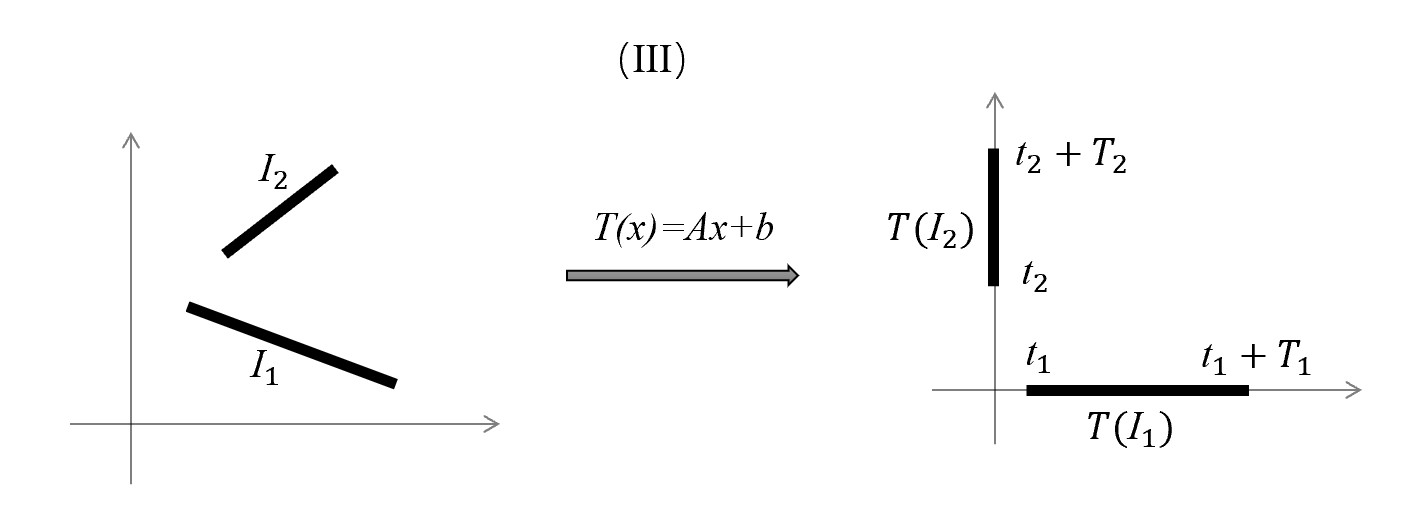} 
  \caption{(III): Affine transformation} \label{twolineself}
\end{figure}

This leads to the following natural question:
\begin{que}
    Given any two line segments in the plane, when is the arc-length measure supported on their union a spectral measure, and, if it is spectral, must all its spectra be line spectra?
\end{que} 
In this paper, we provide a complete answer to this question. To help the reader follow the discussion, we briefly summarize the known results in this area in Remark \ref{T1}. 
 \begin{remark}\label{T1}
Let $I_1, I_2 \subset \mathbb{R}^2$ be two  non-overlapping (but possibly intersecting)  line segments with $|I_1|+|I_2|=2$. Let $\nu$ be arc-length measure on both $I_1$ and $I_2$ normalized (divided by 2) so as to be a probability measure. We record the following three facts. 
\begin{itemize}
\item[(I)]  {\bf Colinear line segments.} If $I_1$ and $I_2$ are collinear line segments, then, by results of Dutkay and Lai \cite{DL2014} and {\L}aba \cite{laba2001}, the  measure $\nu$ is spectral (i.e.\ the set $I_1 \cup I_2$ is a spectral set in $\RR$) if and only if the gap between $I_1$ and $I_2$ is an integer when $|I_1| = |I_2|$, and the gap between $I_1$ and $I_2$ is an even integer when $|I_1| \neq |I_2|$. (Refer to Fig. \ref{twolineco}) 
\item[(II)]   {\bf Parallel line segments.}  If $I_1$ and $I_2$ are parallel line segments (non-collinear)   of equal length, Kolountzakis and Lai \cite[ \S 6]{kolountzakis2025non} show that whenever  the measure $\nu$ is spectral  it always has a line spectrum but also admits a spectrum that is not contained in any line. (Refer to Fig. \ref{twolineco})
\begin{figure}[htbp]
  \centering
\includegraphics[width=1.0\textwidth]{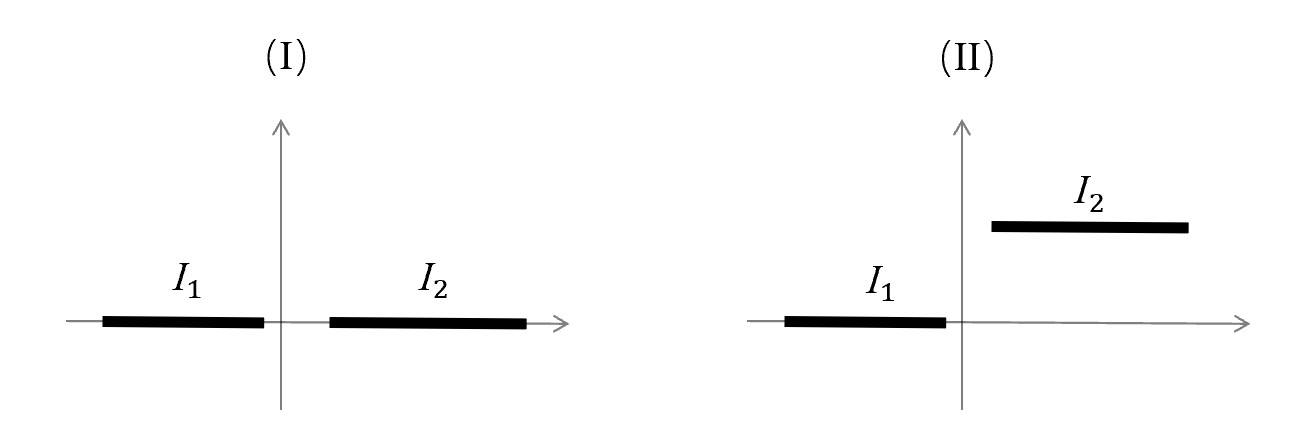} 
  \caption{(I)  and  (II).} \label{twolineco}
\end{figure}
\item[(III)] {\bf Non-parallel line segments.} If $I_1$ and $I_2$ are two non-parallel line segments, then there exists an invertible affine map $T(x)=A x+ b$ with $A \in GL(2, \mathbb{R})$ such that $T\left(I_1\right)$ is contained in the $x$-axis and $T\left(I_2\right)$ is contained in the $y$-axis. Moreover, if $\mu$ is a spectral measure with spectrum $\Lambda \subset \mathbb{R}^2$, then the push-forward measure $\nu=T_{\#} \mu$ is again spectral, with spectrum $(A^{-1})^T \Lambda$. In particular, if $\Lambda$ is contained in a line $L\subset\mathbb{R}^2$, then $(A^{-1})^T\Lambda$ is contained in the line $(A^{-1})^T L$. Thus the property of being spectral is invariant under invertible affine changes of variables, and we may work in this normalized configuration without loss of generality. Notice that, additionally, we can always choose this affine transformation so that $T(I_1)$ and $T(I_2)$ have the same lengths as $I_1$ and $I_2$ (Refer to Fig. \ref{twolineself}). 
In other words, when we study the spectrality of the arc-length measure   supported on two non-parallel line segments in the plane, we may reduce the problem to the study of the spectrality of the   measure
\begin{equation}\label{rho-measure}
\rho=\frac{1}{2}L_{[t_1, t_1+T_1]}\times\delta_0+\frac{1}{2}\delta_0\times L_{[t_2, t_2+T_2]}   \text{ with } T_1+T_2=2.
\end{equation}
\end{itemize}
\end{remark}

According to Remark \ref{T1}, to fully answer the question above, we  still  need to analyze the case of two parallel line segments  (non-collinear) and the case of the   measure $\rho$ defined above  in \eqref{rho-measure}. Our first result is about the spectrality of two parallel line segments. 

\begin{theorem}\label{twoline00} 
Given are real numbers $a_1$, $a_2$  and  $h_1\neq h_2$. Let $T_1, T_2>0$ and $T_1+T_2=2$. Then the measure $\nu=\frac{1}{2}L_{[a_1, a_1+T_1]}\times\delta_{h_1}+\frac{1}{2} L_{[a_2, a_2+T_2]}\times\delta_{h_2}$   is always spectral with one of its spectra contained in a straight line. 
\end{theorem}
 For a case of non-parallel line segments, i.e.  the measure  $\rho$ of the form  \eqref{rho-measure} with 
$T_1=T_2$, Lu \cite{Lu25} proved that $\rho$ is a spectral measure if and only if $t_1+t_2 \in\mathbb{Z} \setminus\{-1\}$ or $t_1-t_2 \in \mathbb{Z} \setminus\{0\}$. However,  in the paper \cite{Lu25},  he  does not identify the possible structures of the line spectra.  In this paper,  we establish  the sufficient and necessary condition of spectrality for general non-parallel line segments and show  the structures of the line spectra:    a spectral measure of non-parallel line segments  has precisely two possible types of line spectra.   Even for the case where $T_1=T_2$, our proof includes a completely different method from Lu's \cite{Lu25}. Now we state our second result. 
\begin{theorem}\label{twoline1}
Let $t_1,t_2\in\mathbb{R}$ and  $T_1,T_2>0$ with $T_1+T_2=2$. Consider the measure
\[
\rho=\frac{1}{2}L_{[t_1, t_1+T_1]}\times\delta_0+\frac{1}{2}\delta_0\times L_{[t_2, t_2+T_2]}.
\]
Then $\rho$ is a spectral measure if and only if one of the following  statements hold.
\begin{itemize}
\item[(I)] When $T_1=T_2$, we have  $t_1+t_2 \in\mathbb{Z} \setminus\{-1\}$ or $t_1-t_2 \in \mathbb{Z} \setminus\{0\}$.
\item[(II)] When $T_1\neq T_2$, we have $t_1+t_2\in 2\mathbb{Z}$ or $t_1-t_2-T_2\in 2\mathbb{Z}$ (equivalently, $t_2-t_1-T_1\in 2\mathbb{Z}$).
\end{itemize}
Moreover, if $\rho$ is a spectral measure, then every spectrum of $\rho$ that contains the origin is contained in either the line $y=x$ or the line $y=-x$. More precisely, if $0\in\Lambda$ is a spectrum of $\rho$, then $\Lambda$ has one of the following forms:
\begin{itemize}
\item[(I)] If $T_1=T_2$, then there exist $\alpha_1, \alpha_2\in(-\tfrac{1}{2}, \tfrac{1}{2}]$ such that
\begin{equation}\label{line-1}
 \Lambda=\begin{cases}
\{(n, -n),(n+\alpha_1, -(n+\alpha_1)): n\in \mathbb{Z}\}, & \text{if } t_1+t_2\in  \mathbb{Z}\setminus\{-1\},\\[6pt]
\{(n,  n),(n+\alpha_2,  n+\alpha_2): n\in \mathbb{Z}\}, & \text{if } t_1-t_2 \in  \mathbb{Z}\setminus\{0\}.
\end{cases}
\end{equation}
\item[(II)] If $T_1\neq T_2$, then
\begin{equation}\label{line-2}
 \Lambda=\begin{cases}
\{(\tfrac{n}{2}, -\tfrac{n}{2}): n\in \mathbb{Z}\}, & \text{if } t_1+t_2\in 2\mathbb{Z},\\[6pt]
\{(\tfrac{n}{2},  \tfrac{n}{2}): n\in \mathbb{Z}\}, & \text{if } t_1-t_2-T_2\in 2\mathbb{Z}.
\end{cases}
\end{equation}
\end{itemize}
\end{theorem}

\begin{remark}\label{T2}
Here we note that $t_1-t_2-T_2\in 2\mathbb{Z}$ is equivalent to $t_2-t_1-T_1\in 2\mathbb{Z}$, since $T_1+T_2=2$.
\end{remark}
 
The observations in Remark~\ref{T1} together with Theorems~\ref{twoline00} and~\ref{twoline1} yield the following result.

\begin{theorem}\label{sum}
Let $\mu$ be the arc-length measure supported on the union of two line segments in the plane. Then $\mu$ is a spectral measure if and only if it admits a line spectrum. Furthermore, in the non-parallel case, the  measure admits only line spectra.
\end{theorem}
 
\begin{remark}\label{rem:projections}
As explained in \cite[\S 4]{kolountzakis2025spectrality} the line spectra of a collection of line segments equipped with the arc-length on them always arise as one-dimensional spectra of the projection of these segments onto a straight line. Thus, in the case of two non-parallel line segments on the straight lines $L_1$ and $L_2$ the possible spectra are on either one of the two angle bisectors of $L_1$ and $L_2$ and therefore it suffices to project the two segments on each of these two lines and check (a) if the projection is injective and (b) that the resulting union of two line segments on the line is spectral. This last condition is known to be equivalent to this projection tiling the line \cite{laba2001}.
\end{remark}

\subsection{No line spectra for some arc-length measures supported on more than two line segments}  
Suppose $\mu$ is a measure on the plane which is arc-length on a finite number of line segments. Such measures have recently been examined for their spectrality \cite{lai2021spectral,kolountzakis2025spectrality,kolountzakis2025non}. In all examples known so far any such spectral measure $\mu$ always had a spectrum contained in a straight line (it could possibly have other spectra too without this property).

It is therefore natural to suspect that this might be the general case. Here we show that this is not so. There are finite collections of line segments in the plane such that their accompanying measure is spectral and each of their spectra is not contained in a straight line.

Notice that for such a measure to have a spectrum contained in a straight line $S$ it is necessary that the orthogonal projection from the support of the measure to $S$ is injective. The reason is that all complex exponentials $e^{2\pi \lambda \cdot x}$ with $\lambda \in S$ are constant along the direction $S^\perp$. If the support of the measure contains two points on the same translate of the line $S^\perp$ then a function on our segments that takes different values on these two points cannot be represented by these exponentials. (See Propositon \ref{propline-2}).

Thus it is rather easy to prove that a collection of line segments does not have a line spectrum: it is enough that on any line $S$ the orthogonal projection is not injective on our set. And if such a collection happens to project injectively onto the line $S$, then for it to have a spectrum contained in $S$ it must also project to a spectral measure on $S$ \cite[Theorem 1.2]{kolountzakis2025spectrality}.

For the existence of the spectrum we will need the following theorem. It describes a situation where the convolution of two spectral measures produces a spectral measure, something which, of course, does not happen in general.

\begin{theorem}\label{th:op}
Suppose $\mu_1$ is a measure on $\RR^d$ and $L \subseteq \RR^k\times\Set{0}^{d-k}$ is a spectrum of $\mu_1$. Suppose also that $\nu_1$ is a measure on $\Set{0}\times\RR^{d-k}$ and $M \subseteq \Set{0}^k\times\RR^{d-k}$ is a spectrum of $\nu_1$.
Then $L+M$ is a spectrum of $\mu_1*\nu_1$.
\end{theorem}

 By using Theorem \ref{th:op} and projection, we have the following theorem. 

\begin{theorem}\label{th:op-2}
There exists a  spectral arc-length measure $\mu$ on $\mathbb{R}^2$, supported on the union of at least three  line segments, such that every spectrum $\Lambda$ of $\mu$ is not contained in a straight  line.
\end{theorem}

\subsection{Size of orthogonal sets of exponentials}
In the previous subsections we  dealt  with very concrete arc-length measures supported on finitely many line segments. For these measures,  every spectrum has strong one dimensional features. For example, in the case of two segments every spectrum is contained in a straight line. In the more general case, with three or more segments, we still need good control on how fast an orthogonal set can grow.  A key underlying ingredient behind these results is always the same: the support is essentially one dimensional, so an orthogonal set cannot be too dense. Informally, the size of any orthogonal set inside a ball of radius $R$ grows at most linearly in $R$.

The purpose of the final part of the paper is to isolate this phenomenon and to state it in a clean, abstract form that applies not only to our concrete examples but to any finite collection  of curves. 
 
\begin{theorem}[Any finite collection of curves]\label{curve1}
 Let $\Gamma$ be a finite union of curves in $\mathbb{R}^d$, each of which is rectifiable and has finite, nonzero length, and let $\rho$ be the arc-length measure on $\Gamma$.   If   $\Lambda \subset \mathbb{R}^d$ is an orthogonal set for $\rho $,
then $\sup _{x \in \mathbb{R}^d}\#(\Lambda\cap B_R(x))= O(R)$ as $R\rightarrow\infty.$
\end{theorem}
 Theorem \ref{curve1}  also yields the following corollary \ref{line segment}.

\begin{corollary}[Any finite collection of line segments]\label{line segment}
Let $d \geq 2$ be an integer.
Let $L_1,\dots,L_m \subset \mathbb{R}^d$ be $m$  line segments, and let $\rho $ be the arc-length  measure on $\bigcup_{k=1}^{m} L_k$. If $\Lambda \subset \mathbb{R}^d$ is an orthogonal set for $\rho$, then $\sup _{x \in \mathbb{R}^d}\#(\Lambda\cap B_R(x))= O(R)$ as $R\rightarrow\infty.$
\end{corollary}


\section{Spectrality of two line segments}\label{s:2-segments}
\subsection{Parallel line segments}
Define the measure $\nu=\frac{1}{2}L_{[a_1, a_1+T_1]}\times\delta_{h_1}+\frac{1}{2} L_{[a_2, a_2+T_2]}\times\delta_{h_2}$   with $h_1\neq h_2$, $T_1, T_2>0$ and $T_1+T_2=2$.  Before proving that $\nu$ always admits line spectra, we first need to prove the following proposition.

\begin{proposition}\label{propline}
Let $I_{1}=[a_1, a_1+T_1]\times\{h_1\}$ 
and  $I_{2}=[a_2, a_2+T_2]\times\{h_2\}$ with $h_1\neq h_2$, $T_1, T_2>0$ and $T_1+T_2=2$. 
For any $k>0$,  there exists a line $L$ (not perpendicular to $I_{1}$ and $I_{2}$) passing through the origin
such that
\begin{itemize}
\item [(a)] the orthogonal projection onto $L$ is injective  on $I_1\cup I_2$, and
\item [(b)] the length of the gap between the orthogonal projections of  $I_1$ and $I_2$  onto the line 
 $L$ is $k$ times the sum of the lengths of their projections.
\end{itemize}
\end{proposition}

\begin{proof}
Since translation of the plane does not change parallelism, lengths, nor the shape of orthogonal projections, we may assume without loss of generality that $a_1=-\frac{T_1}{2}$ and $h_1=0$.
Let $\gamma =(1,0)$, then $$I_1=\left\{t\gamma : t\in [-\frac{T_1}{2}, \frac{T_1}{2}] \right\} \  \ \  \text{and} \ \ \   I_2=\left\{w +t\gamma : t\in  [-\frac{T_2}{2},  \frac{T_2}{2}] \right\}, $$  where $w=(a_2+\frac{T_2}{2}, h_2 )$.
Let $L$ be a line through the origin with unit direction vector $v$ and $\langle \gamma, v \rangle \neq 0$. Then the projections of $I_{1}$ and $I_{2}$ onto the line  $L$ have lengths
$\ell_1 = T_1|\langle \gamma,v\rangle|$ and $\ell_2 = T_2|\langle \gamma,v \rangle|$. Hence the sum of the projected lengths is $\ell_1+\ell_2=2|\langle \gamma,v\rangle|$. (Refer to Fig.\ \ref{proj00}.)
\begin{figure}[htbp]
  \centering
\includegraphics[width=0.8\textwidth]{proj2.jpg} 
  \caption{For every $k>0$ 
 there exists a line 
$L$ through the origin  such that the distance between the projections of 
 $I_1$ and $I_2$  onto the line $L$ is 
 $k$ times the sum of their projected lengths. } \label{proj00}
\end{figure}

For any  $k>0$  we can choose $v_k$ to satisfy the   equation 
$$|\langle w, v_k\rangle|=(2k+1)|\langle \gamma, v_k\rangle|,$$
   that is,  $ v_k$ is perpendicular to $w-(2k+1)\gamma$. So $$|\langle w, v_k\rangle|-  |\langle \gamma, v_k\rangle|= 2k |\langle \gamma, v_k\rangle|,$$
that is, $|\langle w, v_k\rangle|-\frac{1}{2}(\ell_1+\ell_2)= k (\ell_1+\ell_2)$. Note that   $|\langle w, v_k\rangle|-\frac{1}{2}(\ell_1+\ell_2)$ is the length of the gap between the orthogonal projections of  $I_1$ and $I_2$  onto the line $L$.
 
\end{proof}

Let $L$ be a straight line through the origin and $u$ be a unit vector along $L$. Let us also denote by $u^\perp$ the orthogonal subspace to $L$ (a straight line also). We denote by $\pi_L$ the orthogonal projection operator onto line $L$ (but taking values in $\RR$). In other words $\pi_L(v) = t$ for any $v \in t u + u^\perp$.

 Suppose  $\nu$ is a Borel measure on $\RR^2$, then the projection of $\nu$ onto $L$ is the measure $\pi_L\nu$ on $\RR$ defined by
$$
\pi_L\nu(E) = \nu(E u + u^\perp),
$$
where $E \subseteq \RR$.

In \cite[Theorem 1.2]{kolountzakis2025spectrality}, the first and third author of this paper  proved the following very useful theorem, which allows one to detect when a measure on 
$\RR^2$ has a spectrum lying in a straight line by looking only at the spectrum of an appropriate one dimensional projection. We emphasize that this criterion detects only spectra contained in a straight line and cannot be used to disprove the existence of \textit{all} possible spectra.

\begin{theorem}\cite[Theorem 1.2]{kolountzakis2025spectrality}\label{th:projections}
Suppose $\rho$ is a  measure on $\RR^2$ whose support is a finite union of line segments. Suppose also that $L$ is a straight line through the origin such that the orthogonal projection $\pi_L$ onto $L$ is one-to-one $\rho$-almost everywhere. Then $\rho$ has a spectrum $\Lambda u \subseteq L$ if and only if the projection measure $\pi_L\rho$ has spectrum $\Lambda \subseteq \RR$, where $u$ is a unit vector along $L$. 
\end{theorem}
Using Theorem \ref{th:projections} and Proposition \ref{propline} we can now give the proof of Theorem \ref{twoline00}.

\begin{proof}[Proof of Theorem \ref{twoline00}]  
 From Proposition \ref{propline}, for any positive integer $k$, we can find a line $L$ such that the length of the gap between the orthogonal projections of  $[a_1, a_1+T_1]\times\{h_1\}$ 
and  $[a_2, a_2+T_2]\times\{h_2\}$   onto the line 
 $L$ is $k$ times the sum of the lengths of their projections.
 Combining this with Theorem \ref{th:projections} and \cite[Theorem 1.1] {laba2001}, we can deduce the result.

\end{proof}

\begin{remark}\label{rem:thin-spectrum}
The previous proof of Theorem \ref{twoline00} gives an infinity of line spectra for a set of two parallel line segments not on the same line. This infinite collection of different spectra has the additional property that they can be arbitrarily sparse: if the line onto which we are projecting our two segments is almost perpendicular to the intervals themselves then the projection of the two segments has support whose measure tends to 0. Therefore the spectra of the projection, which are also spectra of the initial set of two segments, have linear density that can be arbitrarily small.
\end{remark}
\begin{remark}\label{rem:projection-to-single-interval}
It is worth pointing out that if we only care to show the existence of some spectrum for two parallel line segments not on the same line then we can project them onto a straight line such the projection is injective (with the exception of the endpoints) and the support of the projection is a single interval, which always has a spectrum. To find this straight line suppose without loss of generality that our segments are parallel to the $x$-axis and define the straight line $S$ which joints the left endpoint of one interval with the right endpoint of the other. Project then onto $S^\perp$ to obtain a projection that is just one interval. (Refer to Fig.\ \ref{twopro}.)
\begin{figure}[htbp]
  \centering
\includegraphics[width=0.5\textwidth]{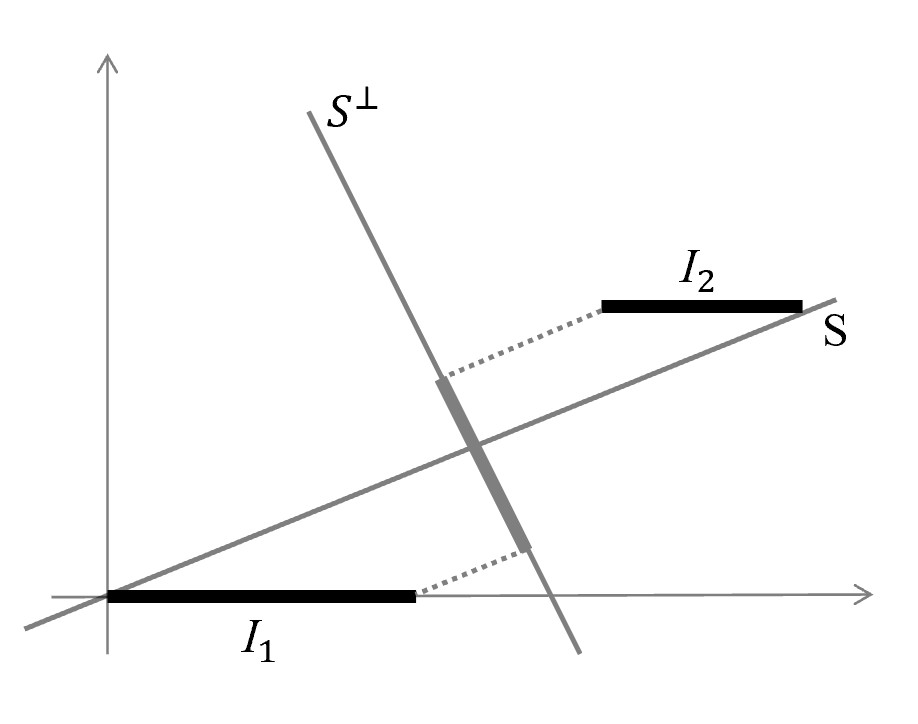} 
  \caption{Projecting to an interval.} \label{twopro}
\end{figure}
\end{remark}

\subsection{Non-parallel line segments}
According to Remark \ref{T1}  (III), when studying the spectrality of two non-parallel line segments we can always reduce to the same question for two line segments that are perpendicular to each other. So let $\rho=\frac{1}{2}L_{[t_1, t_1+T_1]}\times\delta_0+\frac{1}{2}\delta_0\times L_{[t_2, t_2+T_2]}$ with $T_1+T_2=2$ and  $T_1, T_2>0$, and let $\mathcal{Z}\left(\hat{\rho}\right)$ denote the zero set of the Fourier transform $\hat{\rho}$. We know (see for instance \cite{kolountzakis2025spectrality}) that $\left(\lambda_1, \lambda_2\right) \in \mathcal{Z}\left(\hat{\rho}\right)$ if and only if
\begin{align}\label{line1}
e^{\pi i \left(\lambda_1\left(2 t_1+T_1\right)-\lambda_2\left(2 t_2+T_2\right)\right)} \frac{\sin \pi T_1\lambda_1}{\pi \lambda_1}+\frac{\sin \pi T_2\lambda_2}{\pi \lambda_2}=0,
\end{align}
which implies that $ \mathcal{Z}\left(\hat{\rho}\right)=Z_1\cup Z_2$, where
$$
Z_1=   \frac{ \mathbb{Z} \backslash\{0\}}{T_1}\times\frac{ \mathbb{Z} \backslash\{0\}}{T_2}$$
and, defining $T\left(\lambda_1, \lambda_2\right):=\lambda_1\left(2 t_1+T_1\right)-\lambda_2\left(2 t_2+T_2\right)$,
\ifdefined\SMART
\begin{align*}
Z_2 &= \left\{(\lambda_1, \lambda_2 ): T\left(\lambda_1, \lambda_2\right) \in \mathbb{Z}\right.\\
&\left.\hspace{2em}\text{and}\ \  (-1)^{T\left(\lambda_1, \lambda_2\right)} \frac{\sin \pi T_1\lambda_1}{\pi \lambda_1}+\frac{\sin \pi T_2\lambda_2}{\pi \lambda_2}=0\right\}.
\end{align*}
\else
$$
Z_2 = \left\{(\lambda_1, \lambda_2 ): T\left(\lambda_1, \lambda_2\right) \in \mathbb{Z} \ \ \text{and}\ \  (-1)^{T\left(\lambda_1, \lambda_2\right)} \frac{\sin \pi T_1\lambda_1}{\pi \lambda_1}+\frac{\sin \pi T_2\lambda_2}{\pi \lambda_2}=0\right\}.
$$
\fi

 In what follows, we use a sequence of lemmas to show that the projections of the spectrum of $\rho$ onto the $x$-axis and the $y$-axis are both periodic sets. 

\begin{lemma}\label{lemtwo-l--1}
Let $0\in\Lambda$ be a spectrum for $L^2(\rho)$. Then the following two statements hold.
\begin{itemize}
\item[(I)]
  If $T_1\geq T_2$, then  $\Lambda$  has at most one point on each vertical line.
  \item[(II)]  If $T_1\leq T_2$, then  $\Lambda$  has at most one point on each horizontal line.
  \end{itemize}
\end{lemma}
\begin{proof}
We first prove (I). Suppose that $\Lambda $ does not have   at most one point on each vertical line. Then, since the difference of these two points of $\Lambda$ must be a zero of $\ft{\rho}$, it follows from \eqref{line1} that there exist  $(\lambda_1, \lambda_2)\in\Lambda$ and $(\lambda_1, \lambda_2')\in\Lambda$ with $\lambda_2\neq\lambda_2'$ such that
$$ e^{-\pi i(\lambda_2-\lambda_2')\left(2 t_2+T_2\right)  }    T_1 +\frac{\sin \pi T_2(\lambda_2-\lambda_2')}{\pi (\lambda_2-\lambda_2')}=0,$$
which means $\Abs{\frac{\sin \pi T_2(\lambda_2-\lambda_2')}{\pi T_2(\lambda_2-\lambda_2')}}=\frac{T_1 }{T_2}$. Note that $\Abs{\frac{\sin \pi T_2(\lambda_2-\lambda_2')}{\pi T_2(\lambda_2-\lambda_2')}}<1\leq\frac{T_1}{T_2}$,  which is a  contradiction.
For (II), we also arrive at this conclusion through a proof similar to (I).

\end{proof}
 We now show that any spectrum of $L^2(\rho)$ is contained in a family of equidistant parallel lines by the next two lemmas.
 
\begin{lemma}\label{lemtwo-l--2}
Let $0\in\Lambda$ be a spectrum for $L^2(\rho)$. Then $\Lambda \not\subseteq \frac{\mathbb{Z}}{T_1}\times\frac{\mathbb{Z}}{T_2}$ and  
 $$\Lambda\subseteq \left\{(\lambda_1, \lambda_2 ): T( \lambda_1, \lambda_2 ) =\lambda_1\left(2 t_1+T_1\right)-\lambda_2\left(2 t_2+T_2\right)\in \mathbb{Z}\right\}.$$ 
\end{lemma}
\begin{proof}
 Let $$ H_1=\left\{(\lambda_1, \lambda_2 ):  \lambda_1\left(2 t_1+T_1\right)-\lambda_2\left(2 t_2+T_2\right)\in \mathbb{Z}\right\} $$
and
$H_2:=\frac{\mathbb{Z}}{T_1}\times\frac{\mathbb{Z}}{T_2}$.
 $H_1$ and $H_2$ are both subgroups of $\RR^2$. 
 It follows from Lemma 11.4 in \cite{GL17} that $ \Lambda \subseteq H_1$ or $\Lambda \subseteq H_2$.
Suppose that $\Lambda \subseteq H_2$.

If $T_1\geq T_2$, then we 
take the function $f \in L^2(\rho)$ which is $1$ on the horizontal segment and $0$ on the vertical segment.  Notice that it does not matter if the two segments of $\rho$ intersect as they will intersect on at most one point which has $\rho$-measure 0. Then, for any $\lambda=(\lambda_1, \lambda_2 )\in\Lambda$, we have $\lambda_1\in\frac{\mathbb{Z}}{T_1}$ and 
$$
\begin{aligned}
\left\langle f, e_\lambda\right\rangle & =\int f(x, y) e^{-2 \pi i\left(\lambda_1 x+\lambda_2 y\right)} d \rho(x, y) \\
& =\frac{1}{2} \int_{t_1}^{t_1+T_1} e^{-2 \pi i \lambda_1 x} d x \\
& = \begin{cases} \frac{T_1}{2} & \text { if } \lambda_1=0; \\
0 & \text {otherwise.}\end{cases} 
\end{aligned}
$$
Moreover, Lemma \ref{lemtwo-l--1} (I) implies that there is at most one point of $\Lambda$ with $\lambda_1=0$, so this point is the origin, which means $  f=\sum_{\lambda \in \Lambda}\left\langle f, e_\lambda\right\rangle e_\lambda$ is constant $\rho$-almost everywhere  (the series has one term only), a contradiction.  Thus $\Lambda \subset H_1$. 

If $T_1\leq T_2$, then we 
take the function $f \in L^2(\rho)$ which is $1$ on the vertical segment and $0$ on the horizontal  segment  and we can obtain the conclusion  by   a proof similar to the case $T_1\geq T_2$. 

\end{proof}
 
\begin{lemma}\label{lemtwo-l--2-3}
Let $0\in\Lambda$  be a spectrum for $L^2(\rho)$. Then $\Lambda$ is distributed on equidistant parallel lines.
\end{lemma}
\begin{proof}
It follows from Lemma \ref{lemtwo-l--2} that  we have$$\Lambda\subseteq \left\{(\lambda_1, \lambda_2 ): T( \lambda_1, \lambda_2 ) =\lambda_1\left(2 t_1+T_1\right)-\lambda_2\left(2 t_2+T_2\right)\in \mathbb{Z}\right\}.$$  Moreover, 
 \cite[Theorem 1]{kolountzakis2025non} shows  $2 t_1+T_1$ and $2 t_2+T_2$ cannot both be equal to zero. Otherwise, $\rho$ is not a spectral measure. So $\Lambda$ is distributed on equidistant parallel lines. (Refer to Fig. \ref{lines-3}.)
\begin{figure}[htbp]
  \centering
\includegraphics[width=1\textwidth]{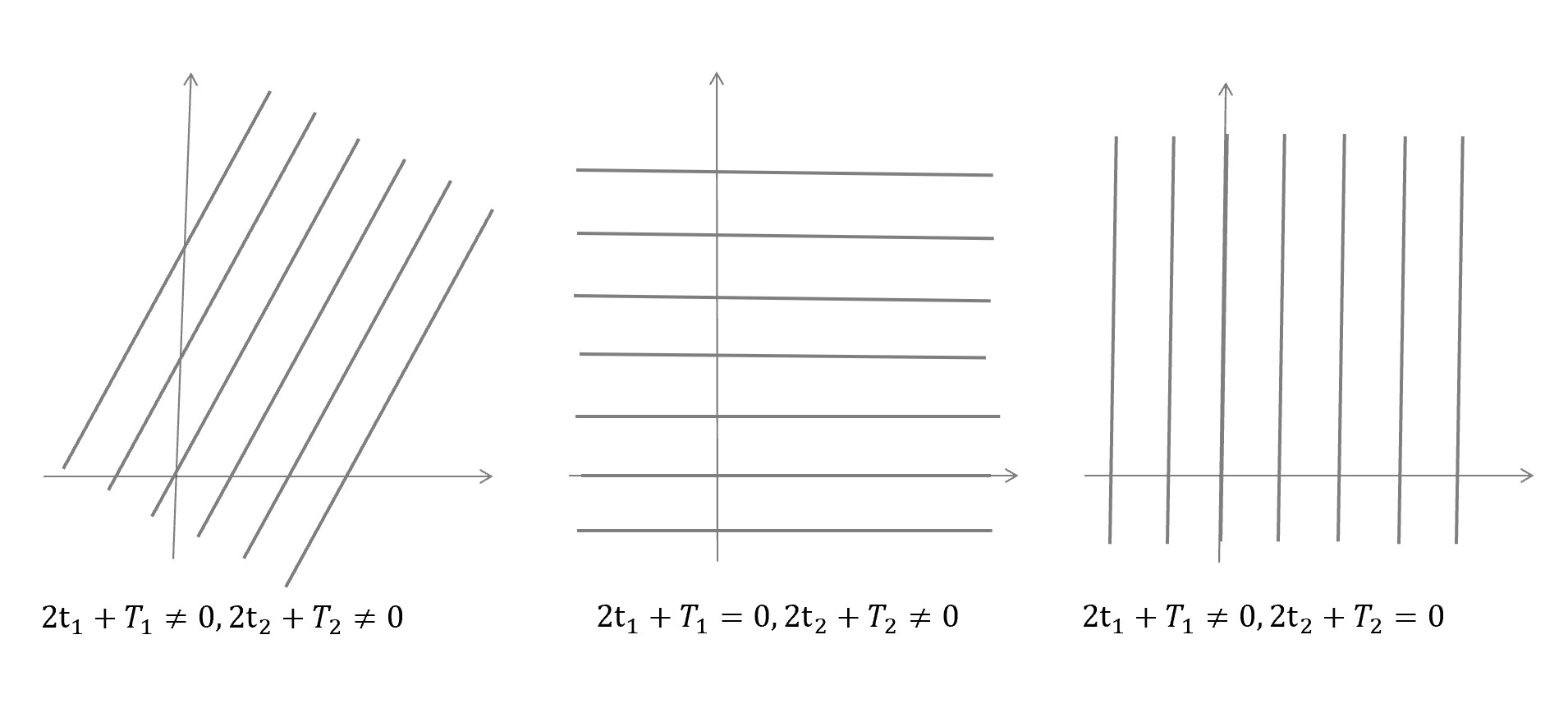} 
  \caption{Three kinds of   equidistant parallel
lines.}  \label{lines-3}
\end{figure}

\end{proof}

We define  the projections  of $\Lambda$ on the $x$-axis and $y$-axis as
$$\Lambda_x=\{\lambda_1: (\lambda_1,\lambda_2 )\in\Lambda\} \ \ \text{and} \  \ \Lambda_y=\{\lambda_2: (\lambda_1,\lambda_2 )\in\Lambda\},$$ respectively. These sets  can be  multi-sets.
\begin{lemma}\label{twoline1-2} If $\Lambda$ is a spectrum of $\rho$,  then  
\begin{align}\label{lem1-2-1}
\frac{2}{T_1}=\sum_{\lambda \in T_1\Lambda_x}\left|\hat{ \mathbf{1}}_{\left[0, 1\right]}\right|^2\left(s-\lambda\right)
\ \ \  \ \text{and} \ \ \ \ \frac{2}{T_2}=\sum_{\lambda \in T_2\Lambda_y}\left|\hat{ \mathbf{1}}_{\left[0, 1\right]}\right|^2\left(s-\lambda\right) 
\end{align} 
for   $s \in \mathbb{R}$.
Additionally,
\begin{enumerate}
\item [(I)] If  $\Lambda$  has at most one point on each vertical line, then $\Lambda_x$ is a periodic set.
\item [(II)] If  $\Lambda$  has at most one point on each horizontal line, then $\Lambda_y$ is a periodic set.
\end{enumerate}

\end{lemma}

\begin{proof}
Assume first that $\Lambda$  has at most one point on each vertical line. Consider the function $f \in L^2(\rho)$ which is $0$ on the $y$-axis and equal to $g(x)=e^{2 \pi i x \frac{s}{T_1}}$ on $[t_1, t_1+T_1] \times\{0\}$ .
 It follows from Parseval's theorem and $f(x, y) \in L^2\left(\rho \right)$ that

$$
\begin{aligned}
\frac{T_1}{2} &= \frac{1}{2} \int_{t_1}^{t_1+T_1}|g(x)|^2 d x\\
   &= \int|f(x, y)|^2 d \rho\\
   & =\sum_{\left(\lambda_1, \lambda_2\right) \in \Lambda}\left|\frac{1}{2} \int_{t_1}^{t_1+T_1} f(x, 0) e^{-2 \pi i \lambda_1 x} d x\right|^2 \\
   & =\frac{1}{4} \sum_{\lambda_{1}  \in \Lambda_x}\left|\int_{t_1}^{t_1+T_1} g(x) e^{-2 \pi i \lambda_1 x} d x\right|^2
\end{aligned}
$$
which gives
\ifdefined\SMART
\begin{align*}
2T_1 &= \sum_{\lambda_{1}  \in  \Lambda_x}\left|\int_{t_1}^{t_1+T_1}e^{2 \pi i x \frac{s}{T_1}} e^{-2 \pi i \lambda_{1}  x} d x\right|^2\\
&= \sum_{\lambda_{ 1}  \in  \Lambda_x}\left|\int_{0}^{T_1}e^{2 \pi i x \frac{s}{T_1}} e^{-2 \pi i \lambda_{1} x} d x\right|^2\\
&= T_1^2\sum_{\lambda_1 \in  \Lambda_x}\left|\int_{0}^{1}e^{2 \pi i x s} e^{-2 \pi i T_1\lambda_1 x} d x\right|^2.
\end{align*}
\else
$$
2T_1 = \sum_{\lambda_{1}  \in  \Lambda_x}\left|\int_{t_1}^{t_1+T_1}e^{2 \pi i x \frac{s}{T_1}} e^{-2 \pi i \lambda_{1}  x} d x\right|^2
= \sum_{\lambda_{ 1}  \in  \Lambda_x}\left|\int_{0}^{T_1}e^{2 \pi i x \frac{s}{T_1}} e^{-2 \pi i \lambda_{1} x} d x\right|^2
= T_1^2\sum_{\lambda_1 \in  \Lambda_x}\left|\int_{0}^{1}e^{2 \pi i x s} e^{-2 \pi i T_1\lambda_1 x} d x\right|^2.
$$
\fi
Therefore,
$$
\frac{2}{T_1}=\sum_{\lambda \in T_1\Lambda_x}\left|\hat{ \mathbf{1}}_{\left[0, 1\right]}\right|^2\left(s-\lambda\right).
$$
By Lemma \ref{lemtwo-l--2}  and Lemma \ref{lemtwo-l--2-3},   we know $\Lambda \not\subseteq \frac{\mathbb{Z}}{T_1}\times\frac{\mathbb{Z}}{T_2}$ and $\Lambda $  is distributed on equidistant parallel lines. Following the arguments in the proofs of Lemmas 5.1, 5.3 and 5.4 in the recent paper \cite{kolountzakis2025spectrality} of the first and third author of this paper, we obtain that there exists a constant  $K>1$ such that for all $\lambda = (\lambda_1, \lambda_2) \in \Lambda \setminus \{(0,0)\}$ we have $K^{-1}\left|\lambda_1\right| \leq\left|\lambda_2\right| \leq K\left|\lambda_{1}\right|$  and $\Lambda_x$ is a set  of finite complexity,  so  $\Lambda_x$ is a  periodic set.

The proof for the case when $\Lambda$ has at most one point on each horizontal line  is essentially identical.

\end{proof}

 \begin{remark}\label{rem:tiling}
The tiling conditions \eqref{lem1-2-1} can also be expressed as follows:
 \begin{align}\label{lem1-2-1-1}
  2 T_1 =\sum_{\lambda \in \Lambda_x}\left|\hat{ \mathbf{1}}_{\left[0, T_1\right]}\right|^2\left(s-\lambda\right)
\ \ \  \ \text{and} \ \ \ \  2 T_2 =\sum_{\lambda \in \Lambda_y}\left|\hat{ \mathbf{1}}_{\left[0, T_2\right]}\right|^2\left(s-\lambda\right).
\end{align}
\end{remark}
 The next lemma provides a convenient normal form for discrete periodic multi-sets. 
\begin{lemma}\label{period}
Let $T, w > 0$. Suppose that $0 \in \Lambda $ is  a discrete periodic  multi-set  and   $$2w=\sum_{\lambda  \in \Lambda}\left|\hat{ \mathbf{1}}_{\left[0, T\right]}\right|^2\left(s-\lambda\right) $$ for $s\in\RR$. Then $\Lambda$  can be represented  as 
 $$\Lambda=\{\lambda^{1}, \lambda^{2},\cdots, \lambda^{m}\}+\frac{mT}{2w}\mathbb{Z} $$
for some  $0 < m\in\ZZ$ and $0 = \lambda^{1}\leq\lambda^{2}\leq\cdots\leq\lambda^{m}<\frac{mT}{2w}$.
\end{lemma}
\begin{proof}
Since $\Lambda $ is  a discrete periodic  multiset  there exists $A>0$ such that $\Lambda=\Lambda+A$   and  the number of points of $ \Lambda$ in any interval of length $ A$ is finite. Hence, there exists a positive integer $m$ such that    $\{\lambda^{1}, \lambda^{2},\cdots, \lambda^{m}\}\subset[0,  A)$ with $0=\lambda^{1}\leq\lambda^{2}\leq\cdots\leq\lambda^{m}<A$. So $ \Lambda=\bigcup_{i=1}^{m}(\lambda^{i}+   A\mathbb{Z})$ and 
\begin{align*}
 2 w &= \frac{1}{ A} \int_0^{A} \sum_{\lambda \in  \Lambda } \left|\hat{ \mathbf{1}}_{\left[0, T\right]}\left(s-\lambda\right)\right|^2 \, ds\\
&=\frac{1}{  A} \sum_{i=1}^{m} \sum_{k\in  \mathbb{Z}} \int_0^{  A} \left|\hat{ \mathbf{1}}_{\left[0, T \right]}\left(s-\lambda^i- Ak\right)\right|^2 \, ds
\\&= \frac{1}{  A} \sum_{i=1}^{m} \sum_{k\in  \mathbb{Z}} \int_{- Ak}^{- A(k-1)} \left|\hat{ \mathbf{1}}_{\left[0, T \right]}\left(s-\lambda^i\right)\right|^2 \, ds
\\&= \frac{1}{  A} \sum_{i=1}^{m} \int_\mathbb{R} \left|\hat{ \mathbf{1}}_{\left[0, T \right]}\left(s-\lambda^i\right)\right|^2 \, ds
\\&= \frac{1}{  A} \sum_{i=1}^{m} \int_\mathbb{R} \left| \mathbf{1}_{\left[0, T\right]}\left(s-\lambda^i\right)\right|^2 ds \ \ \ (\text{Parseval})
\\&= \frac{mT }{ A}.
\end{align*}
This means that $ A= \frac{mT}{2w}$  and $\Lambda=\{\lambda^{1}, \lambda^{2},\cdots, \lambda^{m}\}+\frac{mT}{2w}\mathbb{Z}$. 

\end{proof}
 The next theorem is a simple but powerful rigidity result in one dimension.
It says that if a periodic set $\Lambda$ satisfies the Parseval-type identity coming from the interval $[0,T]$, then $\Lambda$ has essentially only one possible form: it must be the integer translates of two points, and when $T>1$, it must be the half-integer lattice.  
\begin{theorem}\label{twoline1-4}
 Let $T\geq1$ and $m$ be a positive integer. If   $\Lambda =\{\alpha_1, \alpha_2, \cdots ,\alpha_{2m}\}+m\mathbb{Z}$
with   $0=\alpha_1< \alpha_2< \cdots \alpha_{2m}<m $ and  $$
\sum_{\lambda  \in \Lambda }\left|\hat{ \mathbf{1}}_{\left[0, T \right]}\right|^2\left(s-\lambda\right) = 2T $$ for   $s \in \mathbb{R}$, then $ \Lambda=\{0, \alpha\}+\mathbb{Z}$  for some $0<\alpha<1$. Moreover,  we have  $ \Lambda=   \frac{\mathbb{Z}}{2}$ when $T>1$. 
\end{theorem}
\begin{proof}
 Applying the Fourier transform to both sides of the equation 
 $$\sum_{\lambda  \in \Lambda }\left|\hat{ \mathbf{1}}_{\left[0, T \right]}\right|^2\left(s-\lambda\right) = 2T$$
 and the distributional Poisson summation formula, we get
$$
( { \mathbf{1}}_{\left[0, T\right]}* \widetilde{ \mathbf{1}}_{\left[0, T\right]})(x) \cdot \sum_{j=1} e^{2 \pi i \alpha_j x}\cdot \frac{1}{m} \sum_{k \in \mathbb{Z}} \delta_{\frac{k}{m}}=2 T \delta_0.
$$
Since  $\hat{ \mathbf{1}}_{\left[0, T\right]}* \tilde{ \mathbf{1}}_{\left[0, T\right]} > 0$ in $(-T, T)$ and $T\geq 1$, we conclude that
$$
\sum_{j=1}^{2m} e^{2 \pi i \alpha_j\frac{l}{m}}=0\ \ \text{ for any }\ \ l\in\{ \pm 1,  \pm2, \cdots,  \pm  (m-1)\}.
$$  
Let $u_j=e^{2 \pi i \frac{\alpha_j}{m}}$ and  
 denote    by  $p_k\left(u_1, \ldots, u_{2m}\right)$   the $k$-th power sum
$$
p_k\left(u_1, u_2, \ldots, u_{2m}\right)=\sum_{i=1}^{2m} u_i^k=u_1^k+\cdots+u_{2m}^k 
\ \   \text{ for} \ \  k \geqslant 1.
$$
Then we have 
\begin{equation}\label{eq(0.002)}
\begin{cases}
p_1\left(u_1, u_2, \ldots, u_{2m}\right)=0,\\[6pt]
p_2\left(u_1, u_2, \ldots, u_{2m}\right)=0,\\[6pt]
  \cdots  \\
 p_{m-1}\left(u_1, u_2, \ldots, u_{2m}\right)=0. 
 \end{cases}
\end{equation}
The polynomial with roots $u_j$  may be expanded as
$$
P(x):=\prod_{j=1}^{2m}\left(x-u_j\right)=\sum_{j=0}^{2m}(-1)^j e_j x^{2m-j},
$$
where (these are the so-called Newton identities)
\begin{equation*}
\begin{cases} 
e_0 =1, \\[6pt]
e_1 =p_1, \\[6pt]
e_2 =\frac{1}{2}\left(e_1 p_1-p_2\right),\\[6pt]
e_3 =\frac{1}{3}\left(e_2 p_1-e_1 p_2+p_3\right),\\
 \vdots
\end{cases}
\end{equation*}
Here $e_j$ is the $j$-th elementary symmetric function of the numbers $u_1, \ldots, u_{2m}$.

By \eqref{eq(0.002)}, we have $e_1 =e_2 = \cdots =e_{m-1}=0$. This means 
\begin{align*}
P(x) &= x^{2m}+(-1)^{m}e_{m} x^{m}+(-1)^{m+1}e_{m+1} x^{m-1}\\
  &\hspace{3em} +(-1)^{m+2}e_{m+2} x^{m-2}+\cdots+(-1)^{2m}e_{2m}\\
 &= x^{2m}+0\cdot x^{2m-1}+0\cdot x^{2m-2}+\cdots +0\cdot x^{m+1} + (-1)^{m}e_{m}  x^{m} \\
 &\hspace{3em}+(-1)^{m+1}e_{m+1} x^{m-1}+(-1)^{m+2}e_{m+2} x^{m-2}+\cdots+e_{2m}.
\end{align*}
Moreover, $P\left(\frac{1}{x}\right) =\prod_{i=1}^{2m}\left(\frac{1}{x}-u_i\right)=\sum_{j=0}^{2m}(-1)^j e_j x^{j-2m}$ and 
\begin{align*} 
x^{2m}P\left(\frac{1}{x}\right)& = 1+(-1)^{m}e_{m} x^{m} +(-1)^{m+1}e_{m+1} x^{m+1}\\
 &\hspace{3em}+(-1)^{m+2}e_{m+2} x^{m+2}+\cdots+(-1)^{2m}e_{2m}x^{2m}
\\&= e_{2m}x^{2m} - e_{2m-1} x^{2m-2}+\cdots
+ (-1)^{m+1}e_{m+1}x^{m+1}  \\
&\hspace{3em} + (-1)^{m }e_{m } x^{m}  +0\cdot x^{m-1}+\cdots+ 0\cdot x  +1.
\end{align*}
Note that   since the roots of $P(x)$ are $u_1, u_2,\cdots u_{2m}$ and  the roots of $x^{2m}P(\frac{1}{x})$  are $\overline{u_1}, \overline{u_2},\cdots \overline{u_{2m}}$, it follows that  $$P(x)=Cx^{2m}\overline{P}(\frac{1}{x})$$ for $C = e_{2m}$ and  $u_1, u_2,\cdots u_{2m}$  are  the roots of $$x^{2m}+(-1)^{m}e_m x^m+C=0.$$
Since  $u_1=1$ it follows that the numbers $u_1, \ldots, u_{2m}$ are precisely the numbers  $$\{e^{2\pi i  \frac{j}{m}}: j\in\{0,1, \cdots, m-1\}\}\cup\{e^{2\pi i  \frac{\alpha+j}{m}}: j\in\{0,1, \cdots, m-1\}  \}$$ for some $\alpha\in(0, 1)$. It follows that  the numbers $a_1, \ldots, a_{2m}$ are  exactly the numbers
$$0, \alpha, 1, \alpha+1, \ldots, m-1, \alpha+m-1$$ and therefore 
$ \Lambda=\{0, \alpha\}+\mathbb{Z}$.

Moreover, if $T>1$, then 
$$
\sum_{j=1}^{2m} e^{2 \pi i \alpha_j\frac{l}{m}}=0\ \ \text{ for any }\ \ l\in\{ \pm 1,  \pm2, \cdots,  \pm  m\}.
$$
So
\ifdefined\SMART
$$p_1\left(u_1, u_2, \ldots, u_{2m}\right)=\cdots=p_m\left(u_1, u_2, \ldots, u_{2m}\right)=0$$
\else
$$p_1\left(u_1, u_2, \ldots, u_{2m}\right)=p_2\left(u_1, u_2, \ldots, u_{2m}\right)=\cdots=p_m\left(u_1, u_2, \ldots, u_{2m}\right)=0$$
\fi
and $e_1=e_2=\cdots=e_m=0$. From the proof above, we can obtain $ u_1, u_2,\cdots u_{2m}$  are  the roots of $x^{2m}+C=0$. Since  $u_1=1$, this means that $C=-1$ and the numbers $u_1, \ldots, u_{2m}$ are precisely the numbers $$ \{e^{2\pi i\frac{j}{2m}}: j=0, 1, \ldots, 2m-1\}.$$ It follows that  the numbers $a_1, \ldots, a_{2m}$ are the numbers
$0, \frac12, 1, \ldots, m-\frac12$ and therefore $\Lambda = \frac12\ZZ$.

\end{proof}

 The following proposition summarizes the basic  restrictions imposed on $\Lambda-\Lambda$
  by the spectrum condition. These constraints will be used repeatedly in the sequel.
\begin{proposition}\label{lemtw}
Let $0\in\Lambda$  be a spectrum for $L^2(\rho)$,
then the following three  statements hold. 
\begin{itemize}
\item[(I)]
  If $T_1=T_2$, then for any $(\lambda_1, \lambda_2)\in(\Lambda-\Lambda)\setminus\{0\}$ with $|\lambda_1|\leq\frac{1}{2}$, we have $|\lambda_2|=|\lambda_1|$. 
\item[(II)] If $T_1\geq T_2$, then for any $(\lambda_1, \lambda_2)\in(\Lambda-\Lambda)\setminus\{0\}$ with $ |\lambda_1|=\frac{1}{2}$, we have $|\lambda_2|=|\lambda_1|=\frac{1}{2}$. 
\item[(III)] If $T_1\leq T_2$, then for any $(\lambda_1, \lambda_2)\in(\Lambda-\Lambda)\setminus\{0\}$ with $ |\lambda_2|=\frac{1}{2}$, we have $|\lambda_1|=|\lambda_2|=\frac{1}{2}$.
\end{itemize}
\end{proposition}
\begin{proof}
(I). By Lemma \ref{lemtwo-l--2}, we have $T(\lambda_1, \lambda_2 )\in \mathbb{Z}$, which means $
\Abs{\frac{\sin  \pi  T_1\lambda_1 }{ \pi T_1 \lambda_1}} = \Abs{\frac{\sin \pi T_1\lambda_2}{\pi T_1\lambda_2}}.$ The function $f(x)=\frac{\sin   \pi x }{ \pi x } $ shown in the following Figure \ref{fig-1} takes each value only twice when $\Abs{x}\le 1/2$. It follows from $T_1=T_2=1$ that $|\lambda_2|=|\lambda_1|$.
\begin{figure}[htbp]
  \centering
\includegraphics[width=0.7\textwidth]{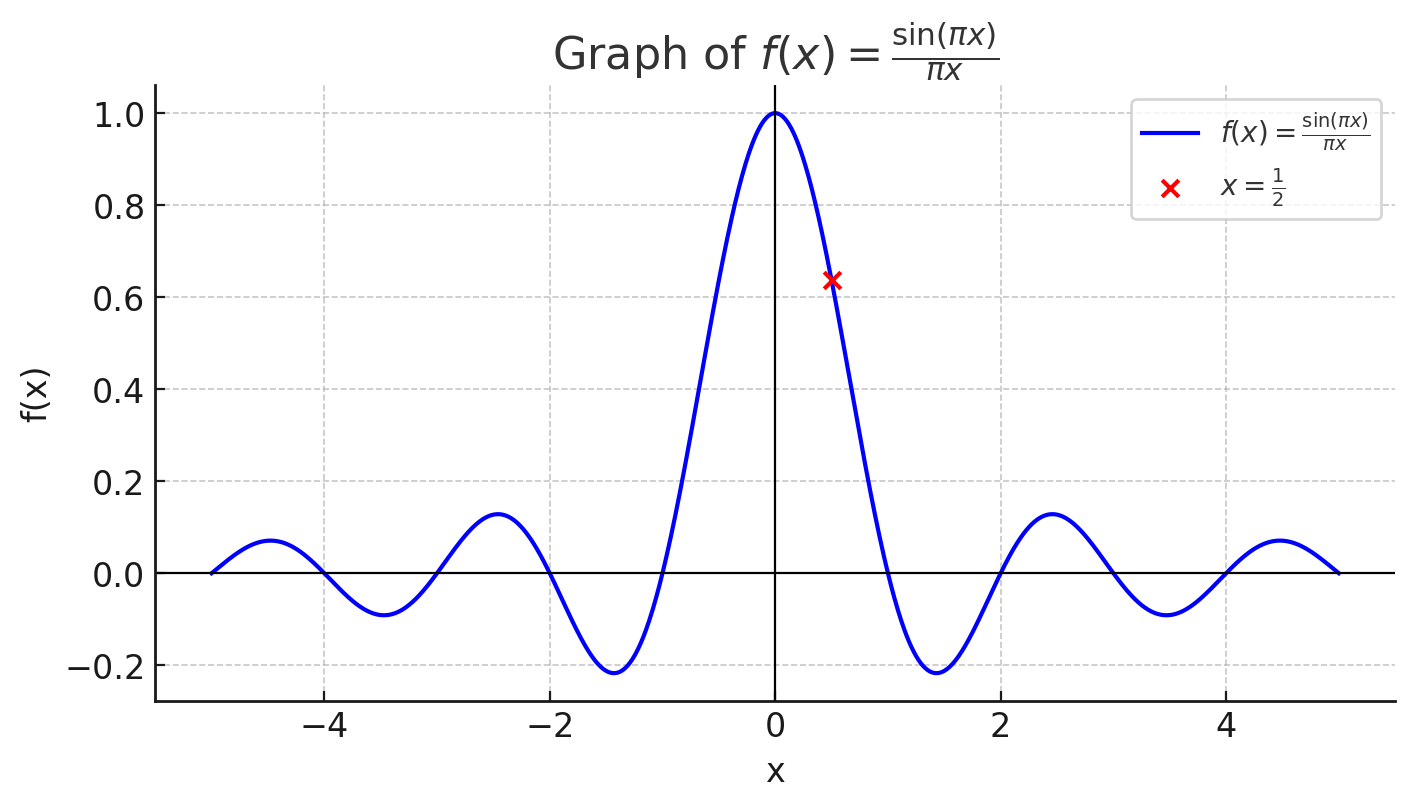} 
  \caption{   $f(x)=\frac{\sin  \pi x}{ \pi x}$ }  \label{fig-1}
\end{figure}
 
(II). If $(\lambda_1, \lambda_2)\in(\Lambda-\Lambda)\setminus\{0\}$ with $ |\lambda_1|=\frac{1}{2}$, by Lemma \ref{lemtwo-l--2}, we have $T(\frac12, \lambda_2 )$  or  $T(-\frac12, \lambda_2 )\in \mathbb{Z}$, and  in any case  we get $
\Abs{\frac{\sin \frac{ 1}{2} \pi  T_1 }{ \frac{ 1}{2}\pi  }} = \Abs{\frac{\sin \pi T_2\lambda_2}{\pi  \lambda_2}}$.
Since $T_1+T_2=2$, we have 
 $\Abs{\frac{\sin \frac{1}{2} \pi  T_2 }{\frac{1}{2}\pi T_2 }} = \Abs{\frac{\sin \pi T_2\lambda_2}{\pi T_2\lambda_2}}.$  Hence, from  $0<\frac{ T_2 }{2}\leq\frac{ 1 }{2} $, we can obtain $|\lambda_2|=|\lambda_1|=\frac{1}{2}$.

(III). Its proof is similar to (II), so we omit it here.

\end{proof}
  With these results in hand, it remains to verify that 
  $\Lambda$  has multiplicity one, i.e.,  $\Lambda$  has at most one point on each vertical line and horizontal line. Once this is established, Lemma \ref{twoline1-2} immediately implies that 
$\Lambda_x$ and $\Lambda_y$ are  periodic sets. In addition, the property of multiplicity one will also play a crucial role later when we determine the precise structure of spectra. 
 
\begin{proposition}\label{lemtwo-l--1-2} 
Let $0\in\Lambda$  be  a spectrum for $L^2(\rho)$. Then \begin{enumerate}
\item [(I)] $\Lambda$  has multiplicity one;
\item [(II)] $\Lambda_x$ and $\Lambda_y$ are  periodic sets.
\end{enumerate}
\end{proposition}
\begin{proof}
(I). If $T_1=T_2$, we can get the conclusion from Lemma \ref{lemtwo-l--1}. If  $T_1>T_2$, then  $\Lambda_x=   \frac{\mathbb{Z}}{2}$ and $\Lambda$ has at most one point on each vertical line  by Theorem \ref{twoline1-4} and Lemma \ref{lemtwo-l--1}.
 Hence, for any $(\lambda_1, \lambda_2)\in\Lambda$ and $(\lambda_1', \lambda_2)\in\Lambda$, we have  $\lambda_1'-\lambda_1'\in\frac{\ZZ}{2}$ and
 $$
\Abs{\frac{\sin (\lambda_1 - \lambda_1')\pi  T_2 }{ (\lambda_1 - \lambda_1')\pi  }}=\Abs{\frac{\sin (\lambda_1 - \lambda_1')\pi  T_1 }{ (\lambda_1 - \lambda_1')\pi  }} = \Abs{\frac{\sin (\lambda_2 - \lambda_2)\pi T_2 }{  (\lambda_2 - \lambda_2)\pi}}=T_2,$$
which means  $ \lambda_1= \lambda_1'$. So $\Lambda$ also has at most one point on each horizontal line. Therefore,  $\Lambda$  has multiplicity one. If  $T_1<T_2$,  we also have that $\Lambda$  has multiplicity one using similar proof methods.

(II). We can get (II) by (I) and Lemma \ref{twoline1-2}.

 \end{proof}

  With multiplicity one and the periodicity of 
 $\Lambda_x$ and 
 $\Lambda_y$	 in hand, we can now pin down the spectrum completely. The next theorem gives the explicit classification of all possible  $\Lambda$.
  
 \begin{theorem} \label{lemtw--t1}
Let $0\in\Lambda$  be a spectrum for $L^2(\rho)$,
then the following two statements hold. 
\begin{itemize}
\item[(I)]
 If $T_1 \neq T_2$, then $ \Lambda= \{(\frac{n}{2}, \frac{n}{2} ): n\in \mathbb{Z}\}$ or $ \Lambda= \{(\frac{n}{2}, -\frac{n}{2} ) : n\in \mathbb{Z}\}$.  
\item[(II)]
 If $T_1=T_2$, then  
$ \Lambda= \{(n, n ), (n+\alpha,n+\alpha): n\in \mathbb{Z}\}$ or $ \Lambda= \{(n, -n ), (n+\alpha, -(n+\alpha)): n\in \mathbb{Z}\}$ for some $\alpha\in(-\frac{1}{2}, \frac{1}{2}] \setminus\{0\}$.  
\end{itemize}
 \end{theorem}
\begin{proof}
(I).  If $T_1> T_2$,  according to Remark \ref{rem:tiling}  and  Theorem \ref{twoline1-4}, we have $\Lambda_x=\{0,  \frac{1}{2}\}+\mathbb{Z}$. And by   Proposition \ref{lemtw}, we have   $(\frac{1}{2}, \frac{1}{2})\in \Lambda$ or  $(\frac{1}{2}, -\frac{1}{2})\in \Lambda$. 
Suppose that $(\frac{1}{2}, \frac{1}{2})\in \Lambda$, then   from Proposition \ref{lemtw} and Proposition \ref{lemtwo-l--1-2}, we have $(-\frac{1}{2}, -\frac{1}{2}),  (1, 1) \ \text{and}\ (- 1, -1)  \in\Lambda$. By induction, assume $(\frac{m}{2}, \frac{m}{2})\in \Lambda$ for all $|m|\leq n$.  If  $|m|=n+1$, there exists  $\gamma$ such that $ (\frac{m}{2}, \gamma)\in\Lambda$. Since $ (\frac{m-1}{2}, \frac{m-1}{2})\in\Lambda$, from Proposition \ref{lemtwo-l--1-2}, we have $ \gamma=\frac{m}{2}$. Hence, $ \Lambda= \{(\frac{n}{2}, \frac{n}{2} ): n\in \mathbb{Z}\}$. Similarly, for the case where $(\frac{1}{2}, -\frac{1}{2})\in \Lambda$, we can get  $ \Lambda= \{(\frac{n}{2}, -\frac{n}{2} ) : n\in \mathbb{Z}\}$.

If  $T_1< T_2$, we can obtain  similarly  that 
 $ \Lambda= \{(\frac{n}{2}, \frac{n}{2} ): n\in \mathbb{Z}\}$ or $ \Lambda= \{(\frac{n}{2}, -\frac{n}{2} ) : n\in \mathbb{Z}\}$.
 
 (II). According to Remark \ref{rem:tiling} and  Theorem \ref{twoline1-4}, $\Lambda_x$ and $\Lambda_y $ can be represented  as  $\Lambda_x  =\{0, \alpha_1\}+\mathbb{Z}$  and $\Lambda_y  =\{0, \alpha_2\}+\mathbb{Z}$  for some $ \alpha_1, \alpha_2\in(-\frac{1}{2}, \frac{1}{2}]\setminus\{0\}$.

Suppose that $\alpha_1=\frac{1}{2}$ or $\alpha_2=\frac{1}{2}$, then from (I), we have $ \Lambda= \{(\frac{n}{2}, \frac{n}{2} ): n\in \mathbb{Z}\}$ or $ \Lambda= \{(\frac{n}{2}, -\frac{n}{2} ) : n\in \mathbb{Z}\}$. 

Suppose that $\alpha_1\neq\frac{1}{2}$ and $\alpha_2\neq\frac{1}{2}$. From Proposition \ref{lemtw}, we have $\alpha_1=\alpha_2$ or $\alpha_1=-\alpha_2$, and for any $(n, \lambda_2), (n+\alpha_1, \lambda_2')\in\Lambda$ with   $n\in\ZZ$, we have $\lambda_2=m$ for some $m\in\ZZ$ and $|\lambda_2'-\lambda_2|=\alpha_1$, which means $\lambda_2'=m+\alpha_1$ or $\lambda_2'=m-\alpha_1$. Moreover, 
 $$
 \Abs{\frac{\sin \pi  \alpha_1  }{  \pi  (n+\alpha_1)}} =  \Abs{\frac{\sin \pi (n+\alpha_1) }{  \pi  (n+\alpha_1)}} = \Abs{\frac{\sin \pi \lambda_2'  }{ \pi\lambda_2'}} = \Abs{\frac{\sin \pi  \alpha_1  }{ \pi\lambda_2'}},$$
so $\lambda_2'=n+\alpha_1$ or $\lambda_2'=-n-\alpha_1$, which means $\lambda_2=n$ or $\lambda_2=-n$. Hence, $|\lambda_2|=|n|$.

If $\alpha_1=\alpha_2$, then for any $(n, \lambda_2)\in\Lambda$, we have 
\begin{equation}\label{llq}\frac{\sin \pi (n-\alpha_1 )}{  \pi    (n-\alpha_1 ) }  =  \frac{\sin \pi (\lambda_2-\alpha_1 ) }{  \pi  (\lambda_2-\alpha_1)}\ \ \text{or} \ \ \frac{\sin \pi (n-\alpha_1 )}{  \pi    (n-\alpha_1 ) }  = - \frac{\sin \pi (\lambda_2-\alpha_1 ) }{  \pi  (\lambda_2-\alpha_1)}.\end{equation}
 Combining \eqref{llq} with $|\lambda_2|=|n|$, we have $\lambda_2=n$.
These imply that $\{(n, n): n\in\ZZ\}\cup\{(\alpha_1, \alpha_1)\}\subset\Lambda$. So  for any $n\in\ZZ $, we have $(n, n), (n+\alpha_1 , \beta) \in\Lambda$ for some  $\beta\in\RR$. From Proposition \ref{lemtw}, we have $ \beta= n+\alpha_1$  or $\beta=n-\alpha_1$. Suppose that $\beta=n-\alpha_1$, we have 
$$
 \Abs{\frac{\sin n\pi  }{  n\pi   }} =  \Abs{\frac{\sin (n-2\alpha_1)\pi  }{ (n-2\alpha_1) \pi  }},$$
 which contradicts $\alpha_1\notin\frac{\ZZ}{2}$. Therefore, $ \Lambda= \{(n, n ), (n+\alpha,n+\alpha): n\in \mathbb{Z}\}$, where $\alpha=\alpha_1$.

If  $\alpha_1\neq\alpha_2$,  similar to $\alpha_1=\alpha_2$,  we can obtain  $ \Lambda= \{(n, -n ), (n+\alpha, -(n+\alpha)): n\in \mathbb{Z}\}$.

Combining the above, we have completed the  proof.

\end{proof} 
 The occurrence of special points such as 
  $(\frac{1}{2}, \pm\frac{1}{2})$,   $(1, \pm 1)$ in $\Lambda$
  imposes rigid conditions on $t_1$ and $t_2$, and therefore controls the relative placement of the two segments. The following two lemmas record  these. 
\begin{lemma}\label{lemtwo}
Suppose that $\rho $ is a spectral measure with the spectrum $0\in\Lambda$, then the following two statements hold.
\begin{itemize}
\item[(I)] 
If $(\frac{1}{2}, -\frac{1}{2})\in\Lambda$, then $t_1+t_2\in 2\mathbb{Z}$.
\item[(II)]
 If $(\frac{1}{2}, \frac{1}{2})\in\Lambda$, then $t_1-t_2-T_2\in 2\mathbb{Z}$.
 \end{itemize}
\end{lemma}
\begin{proof} First, we prove that $T (\frac{1}{2}, -\frac{1}{2}) , T(\frac{1}{2},  \frac{1}{2}) \in 2\mathbb{Z}+1$. Suppose that $T (\frac{1}{2}, -\frac{1}{2})\in 2\mathbb{Z}$ or $T(\frac{1}{2},  \frac{1}{2}) \in 2\mathbb{Z}$. Then
$$ \sin \frac{\pi T_1}{2}= - \sin \frac{\pi T_2}{2}  = - \sin \frac{\pi (2-T_1)}{2}=-\sin \frac{\pi T_1}{2},$$
which is a contradiction.

  (I). If $(\frac{1}{2}, -\frac{1}{2})\in\Lambda$,  we have $$T (\frac{1}{2}, -\frac{1}{2})= \frac{1}{2}\left(2 t_1+T_1\right)+\frac{1}{2}\left(2 t_2+T_2\right)\in 2\mathbb{Z}+1.$$
Since $T_1+T_2=2$, we have $ t_1+   t_2 \in 2\mathbb{Z}$.

  (II). If $(\frac{1}{2},  \frac{1}{2})\in\Lambda$,  we have $$T (\frac{1}{2}, \frac{1}{2})= \frac{1}{2}\left(2 t_1+T_1\right)-\frac{1}{2}\left(2 t_2+T_2\right)\in 2\mathbb{Z}+1.$$
So $ t_1- t_2- T_2 \in 2\mathbb{Z}$.

\end{proof}
 
\begin{lemma}\label{lemtwo1-5}
Suppose that $\rho $ is a spectral measure with the spectrum $0\in\Lambda$ and $T_1=T_2$. Then the following two statements hold.
\begin{itemize}
\item[(I)] If $(1, -1)\in\Lambda$, then $t_1+t_2\in \mathbb{Z}$. 
\item[(II)]  If $(1, 1)\in\Lambda$, then $t_1-t_2 \in  \mathbb{Z}$.
\end{itemize}
\end{lemma}
\begin{proof}
(I). By Theorem \ref{lemtw--t1}, we have $(\alpha, -\alpha),(\alpha-1, -(\alpha-1)) \in\Lambda$  for some $\alpha$. By \eqref{line1}, we can get $ T(\alpha, -\alpha) , (\alpha-1, -(\alpha-1)) \in 2\mathbb{Z}+1$. Hence, 
$$T (1, -1) =T(\alpha, -\alpha)-T(\alpha-1, -(\alpha-1)) \in 2\mathbb{Z},$$ 
which implies $ t_1+t_2\in\ZZ$.
Similarly, we can get (II).

\end{proof}

 Having obtained the integrality conditions from Lemma 7, we next identify symmetric choices of 
$t_1$ and $t_2$
that are incompatible with spectrality.
 
\begin{lemma}\label{lemtwo1-5-1}
Suppose that  $T_1=T_2$ and one of the following three statements hold. \begin{itemize}
\item[(I)]  $t_1=t_2=-\frac{1}{2}$; 
\item[(II)]  $t_1-t_2=0$ and $t_1+t_2\notin \ZZ$; 
\item[(III)] 
  $t_1+t_2=-1$ and $t_1-t_2\notin \ZZ$.
\end{itemize}
Then  $\rho $ is not a spectral measure.
\end{lemma}
\begin{proof}  In fact, in (I), (II) and (III), the measure $\rho $ is symmetric. See Figure \ref{123pictures}. The conclusion can be directly derived from \cite[Theorem 1]{kolountzakis2025non} and   \cite{kolountzakis2025spectrality}.
\begin{figure}[htbp]
  \centering
\includegraphics[width=1\textwidth]{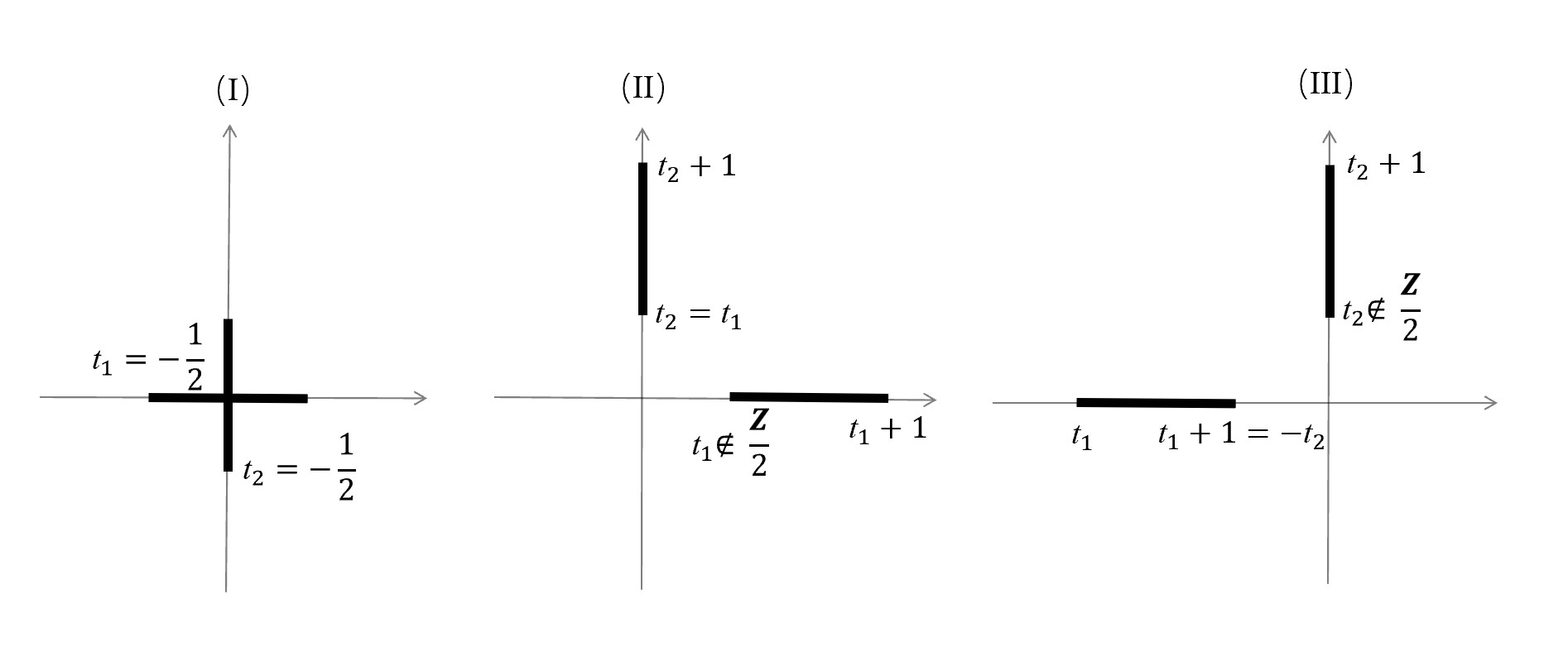} 
  \caption{(I), (II) and (III)  are symmetric measures.} \label{123pictures}
\end{figure}

\end{proof}

 Summarizing the four results above, we can derive the following theorem. 
 \begin{theorem} \label{lemtw--t1-2}
Let $0\in\Lambda$  be a spectrum for $L^2(\rho)$.
Then the following two statements hold. 
\begin{itemize}
\item[(I)] 
  If $T_1=T_2$, then $ t_1+t_2 \in\mathbb{Z} \backslash\{-1\}$ or $  t_1-t_2 \in \mathbb{Z} \backslash\{0\}.$      
\item[(II)]  If $T_1\neq T_2$, then $t_1+t_2\in 2\mathbb{Z}$ or $t_1-t_2-T_2\in 2\mathbb{Z}$. \end{itemize} 
 \end{theorem}
 
\begin{proof}
 Combining Theorem \ref{lemtw--t1}, Lemma \ref{lemtwo}, Lemma \ref{lemtwo1-5} and Lemma\ref{lemtwo1-5-1},
we can obtain this conclusion.

\end{proof}
  To apply \L aba's result \cite[Theorem 1.1] {laba2001}, we will work with the orthogonal projections of $\mathrm{supp}\,\rho$ onto the lines $y=\pm x$. 
We therefore need a  description of the distance (gap) between the two projected components, and it is easy to get the following lemma.
 
\begin{lemma} \label{pro-2}
 Let $I_1= [t_1, t_1+T_1]\times\delta_0$ and $ I_2=\delta_0\times [t_2, t_2+T_2]$ with $T_1, T_2>0$ and $T_1+T_2=2$, then the following four statements hold.
 \begin{itemize}
\item[(I)] 
 If $T_1=T_2$ and $ t_1+t_2 \in\mathbb{Z} \backslash\{-1\}$,  then the length of the gap between the orthogonal projections of  $I_1$ and $I_2$  onto the line 
 $y=-x$ is $\frac{1}{2}(|t_1+t_2+1|-1)$ times the sum of the lengths of their projections. 
\item[(II)]  If $T_1=T_2$ and  $  t_1-t_2 \in \mathbb{Z} \backslash\{0\},$  then the length of the gap between the orthogonal projections of  $I_1$ and $I_2$  onto the line 
 $y=x$ is $\frac{1}{2}(|t_1-t_2|-1)$ times the sum of the lengths of their projections.   
\item[(III)]   If $T_1\neq T_2$ and $t_1+t_2\in 2\mathbb{Z}$, then the length of the gap between the orthogonal projections of  $I_1$ and $I_2$  onto the line 
 $y=-x$ is $\frac{1}{2}(|t_1+t_2+1|-1)$ times the sum of the lengths of their projections.
\item[(IV)] . If $T_1\neq T_2$ and $t_1-t_2-T_2\in 2\mathbb{Z}$, then the length of the gap between the orthogonal projections of  $I_1$ and $I_2$  onto the line 
 $y=x$ is $\frac{1}{2}(|t_1-t_2-T_2+1|-1)$ times the sum of the lengths of their projections.  \end{itemize} 
 \end{lemma}
\begin{proof}These conclusions can be obtained through some simple calculations, which we will omit here.

\end{proof}
 We now conclude the proof of Theorem \ref{twoline1} by combining all  results established above. 
\begin{proof}[Proof of Theorem \ref{twoline1}]  
The necessity part of the theorem can  be deduced from Theorem  \ref{lemtw--t1-2}.  Combining Lemma \ref{pro-2}, Theorem \ref{th:projections} and \cite[Theorem 1.1] {laba2001}, we can  obtain the sufficiency part. Moreover, if $0\in\Lambda$ is a spectrum, then from Theorem \ref{lemtw--t1}, $\Lambda$ must be a straight line spectrum, and the spectrum must be contained in the line $y=x$ or $y=-x$  as   described in  \eqref{line-1} or \eqref{line-2}.

\end{proof}

\begin{remark}
This shows that the spectra of $\rho$  can only be contained in either the line  $y=x$ or $y=-x$.  Indeed, it is clear that the only orthogonal projections of $\rho$ onto a straight line that produce a measure which is a constant multiple of (one-dimensional) Lebesgue measure on its support (as is necessary for spectrality \cite{DL2014}) are the projections onto the straight lines $y=\pm x$. See Figure \ref{twolines}. The content of this section of the paper is therefore to rule out all other possible non-line spectra. 
\begin{figure}[htbp]
  \centering
\includegraphics[width=0.8\textwidth]{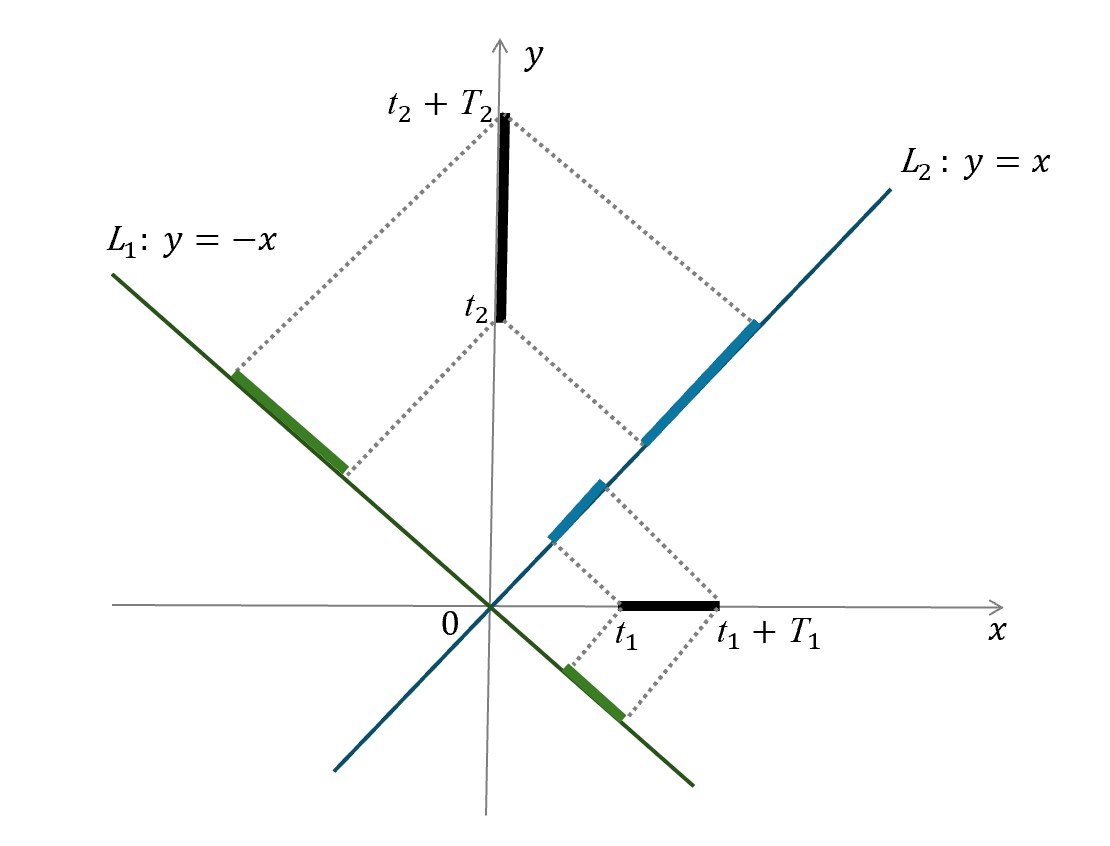} 
  \caption{ There are only two projection methods: $y=x$ or $y=-x$.} \label{twolines}
\end{figure}  
\end{remark}


\section{Spectral collections of segments without line spectra}\label{s:no-line-spectra}

We now turn to the construction of spectral collections of line segments without line spectra.
Theorem~\ref{th:op} plays a crucial role in our subsequent propositions.

For the reader's convenience, we first restate Theorem \ref{th:op} here.

\begin{reptheorem}Suppose $\mu_1$ is a measure on $\RR^d$ and $L \subseteq \RR^k\times\Set{0}^{d-k}$ is a spectrum of $\mu_1$. Suppose also that $\nu_1$ is a measure on $\Set{0}\times\RR^{d-k}$ and $M \subseteq \Set{0}^k\times\RR^{d-k}$ is a spectrum of $\nu_1$.
Then $L+M$ is a spectrum of $\mu_1*\nu_1$.
\end{reptheorem}
\begin{proof} 
Notice that the sums in $L+M$ are all distinct by our assumptions and that we may assume that $\mu_1, \nu_1$ are probability measures. We will verify that
$$
\sum_{\ell\in L, m \in M} \Abs{\ft{\mu_1*\nu_1}}^2(x-\ell-m) = 1
$$
for almost all $x \in \RR^d$.

We have, writing $x=(x_1, x_2)$ with $x_1 \in \RR^k, x_2 \in \RR^{d-k}$,
\ifdefined\SMART
\begin{align*}
\sum_{\ell\in L, m \in M} &\Abs{\ft{\mu_1*\nu_1}}^2(x-\ell-m)\\
  &= \sum_{\ell\in L, m \in M} \Abs{\ft{\mu_1}}^2(x-\ell-m)  \Abs{\ft{\nu_1}}^2(x-\ell-m)\\
 &= \sum_{\ell\in L, m \in M} \Abs{\ft{\mu_1}}^2(x_1-\ell, x_2-m)
    \, \Abs{\ft{\nu_1}}^2(x_1-\ell, x_2-m)\\
 &= \sum_{\ell\in L, m \in M} \Abs{\ft{\mu_1}}^2(x_1-\ell, x_2-m)
    \, \Abs{\ft{\nu_1}}^2(0, x_2-m)\\
 &\vspace{5em}\text{(since $\ft{\nu_1}$ does not depend on the first $k$ coordinates)}\\
 &= \sum_{m \in M}  \Abs{\ft{\nu_1}}^2(0, x_2-m) 
      \sum_{\ell \in L} \Abs{\ft{\mu_1}}^2(x_1-\ell, x_2-m) \\
 &= \sum_{m \in M}  \Abs{\ft{\nu_1}}^2(0, x_2-m) \cdot 1\\
 &= 1.
\end{align*}
\else
\begin{align*}
\sum_{\ell\in L, m \in M} \Abs{\ft{\mu_1*\nu_1}}^2(x-\ell-m)
  &= \sum_{\ell\in L, m \in M} \Abs{\ft{\mu_1}}^2(x-\ell-m)  \Abs{\ft{\nu_1}}^2(x-\ell-m)\\
 &= \sum_{\ell\in L, m \in M} \Abs{\ft{\mu_1}}^2(x_1-\ell, x_2-m)
    \, \Abs{\ft{\nu_1}}^2(x_1-\ell, x_2-m)\\
 &= \sum_{\ell\in L, m \in M} \Abs{\ft{\mu_1}}^2(x_1-\ell, x_2-m)
    \, \Abs{\ft{\nu_1}}^2(0, x_2-m)\\
 &\vspace{5em}\text{(since $\ft{\nu_1}$ does not depend on the first $k$ coordinates)}\\
 &= \sum_{m \in M}  \Abs{\ft{\nu_1}}^2(0, x_2-m) 
      \sum_{\ell \in L} \Abs{\ft{\mu_1}}^2(x_1-\ell, x_2-m) \\
 &= \sum_{m \in M}  \Abs{\ft{\nu_1}}^2(0, x_2-m) \cdot 1\\
 &= 1.
\end{align*}
\fi

\end{proof}

The following Proposition \ref{propline-2} can be found in the  paper  \cite[{\S}4]{kolountzakis2025spectrality} by the first and third author. We include it here only to record a clean formulation.

\begin{proposition}\label{propline-2}
Let $\mu$ be the arc-length measure supported on the union of finitely many line segments
$I_1,\dots, I_N\subset \mathbb{R}^2$. 
Let $L$ be a straight line through the origin and let $\pi_L$ be the orthogonal projection onto $L$.
Assume that $\mu$ is a spectral measure and that it has a spectrum $\Lambda\subset L$.
Then $\pi_L$ is injective $\mu$-almost everywhere on $\supp\mu$.
\end{proposition}
\begin{proof}
Assume that $\pi_L$ is not injective $\mu$-almost everywhere on $\supp\mu$,
 then there exist two disjoint measurable subsets 
$E, F \subset \supp\mu$  such that
$\pi_L(E)=\pi_L(F)=J$. We have 
$$
\pi_{L}\mu|_E=c_j\,H^1|_J\qquad \text{and}
\qquad
\pi_{L}\mu|_F=c_k\,H^1|_J
$$
for some constants $c_j,c_k>0$.
Define
\[
f=c_k\,\mathbf{1}_E-c_j\,\mathbf{1}_F.
\]
Then $f\neq 0$ in $L^2(\mu)$.
For any $\lambda\in\Lambda\subset L$, the exponential
 $ 
e_\lambda(x)=e^{2\pi i\langle \lambda,x\rangle}
 $ 
depends only on $\pi_L x$. Therefore,
\begin{align*}
\langle f,e_\lambda\rangle_{L^2(\mu)}
&=
c_k\int_E e^{-2\pi i\langle \lambda,x\rangle}\,d\mu(x)
-
c_j\int_F e^{-2\pi i\langle \lambda,x\rangle}\,d\mu(x)\\
&=
\int_J e^{-2\pi i\langle \lambda,t\rangle}\,
d(c_k\,\pi_{L}\mu|_E-c_j\,\pi_{L}\mu|_F)(t)\\
&=0.
\end{align*}
Thus $f$ is orthogonal to $\{e_\lambda:\lambda\in\Lambda\}$, contradicting the completeness
of $\Lambda$ in $L^2(\mu)$. Hence $\pi_L$ must be injective $\mu$-almost everywhere on $\supp\mu$.

\end{proof}

We can now use Theorem \ref{th:op} and   Propositon \ref{propline-2}  to show that the following two measures are spectral. One of these two measures consists of parallel line segments and the other contains line segments of different directions.

\begin{proposition}\label{th:parallel}
Suppose $I$ is the line segment from $(0, 0)$ to $(100, 0)$.
Suppose $\mu$ is arc-length on the line segment $S_1 = I$, on $S_2 = I+(0, 1)$ and also on the the line segment $S_3 = I+(\alpha, 2)$, where $\alpha \in (0, \frac{1}{100})$ is an irrational number. Then $\mu$ is spectral but it does not have a line spectrum.
\end{proposition}

\begin{proof}
See Figure \ref{threelines}.

\begin{figure}[htbp]
  \centering
\includegraphics[width=0.7\textwidth]{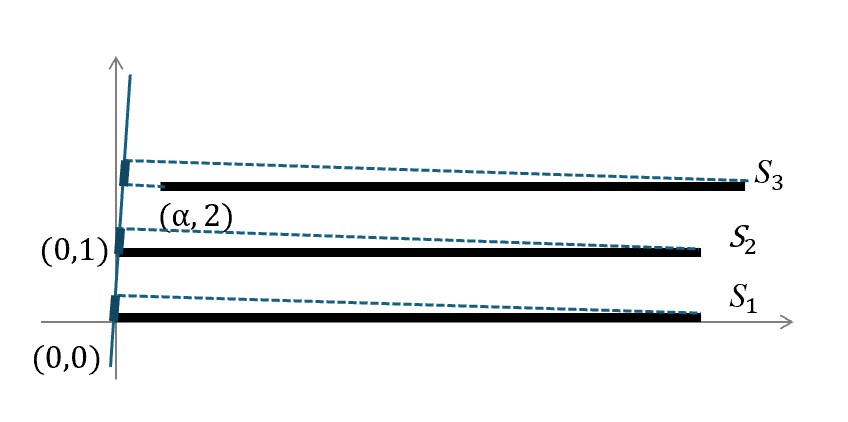} 
  \caption{For irrational number $\alpha$, the arc-length  $\mu$ on the three intervals is  spectral but it does not have a line spectrum. } \label{threelines}
\end{figure}

To prove that $\mu$ is spectral we apply Theorem \ref{th:op} with $\nu_1$ being arc-length on $I$ and $\mu_1 = \delta_{(0, 0)} + \delta_{(0, 1)} + \delta_{(\alpha, 2)}$. Notice that $\mu_1$ has the set $\Set{(0, 0), (0, 1/3), (0, 2/3)}$ as a spectrum.

To show that it does not have a line spectrum, first notice that for $\supp\mu$ to project injectively on the line $S$, the line $S$ must be almost vertical. If this happens it projects to three intervals of equal length, but the two gaps among them have incommensurable lengths because $\alpha$ is irrational. But if this union of three intervals of $S$ is to be spectral each gap they define must be a multiple of the common length of the intervals \cite{ducasse2025spectral,kolountzakis2025geometric}, therefore the projection measure is never spectral. 

\end{proof}

\begin{proposition}\label{th:L}
Suppose $\mu$ is arc-length on the four intervals
\begin{align*}
(0, 0) &- (1, 0),\\
(0, 0) &- (0, 1),\\
(1, 1) &- (1, 2),\\
(1, 1) &- (2, 1).
\end{align*}
Then $\mu$ is spectral but it does not have a line spectrum.
\end{proposition}

\begin{proof}

See Figure \ref{fourlines}.

\begin{figure}[htbp]
  \centering
\includegraphics[width=0.6\textwidth]{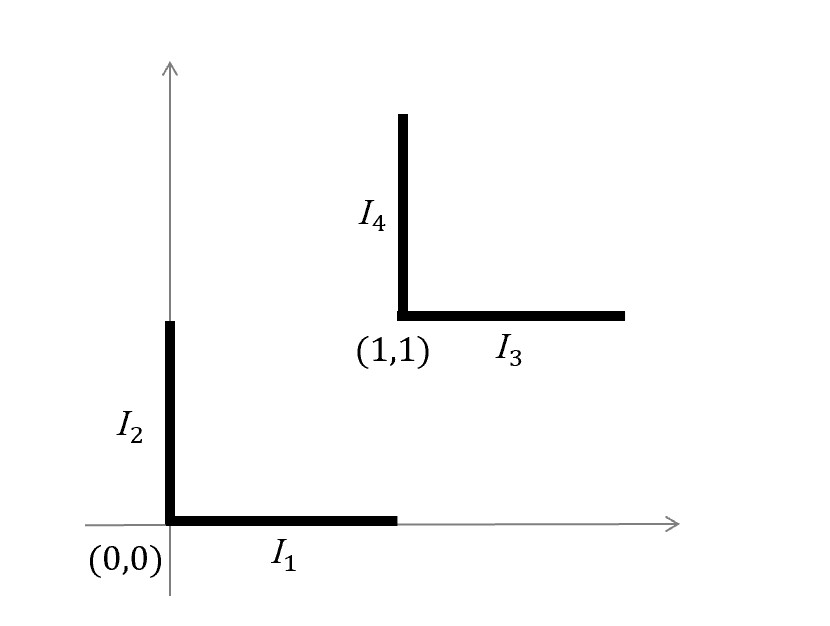} 
  \caption{The arc-length  $\mu$ on the four intervals  is   spectral but it does not have a line spectrum. } \label{fourlines}
\end{figure}

To prove that $\mu$ is spectral we use Theorem \ref{th:op} with $\mu_1$ being the arc-length on the first two line segments above and $\nu_1 = \delta_{(0, 0)} + \delta_{(1, 1)}$. (Here we remember \cite{lai2021spectral} that  $\mu_1$ has a spectrum contained in the line $y=-x$.)

To show that $\mu$ does not have a line spectrum we observe that there is no straight line $S$ onto which $\supp\mu$ projects injectively. 

\end{proof}
\begin{proof}[Proof of Theorem \ref{th:op-2}] The result follows from 
Propositions  \ref{th:parallel} and   \ref{th:L}.

\end{proof}
\section{Growth bound for orthogonal sets in balls}\label{s:size}
 
\subsection{Any finite collection of  curves }
We will use  Lev's a result  \cite[Lemma 3.1.]{Lev18}  to treat the case of arc-length Lebesgue  measures supported on curves. 

 \textbf{Lev's condition  \cite[Lemma 3.1.]{Lev18}}:  Let $\mu$ be a positive, finite measure on $\mathbb{R}^d$. Given a real number $\alpha$ with $0 \leqslant \alpha \leqslant d$, Lev considered  the following condition
 \begin{equation}\label{lev-0}
 \liminf _{R \rightarrow \infty} \frac{1}{R^{d-\alpha}} \int_{|t|<R}|\widehat{\mu}(t)|^2 d t>0
 \end{equation}
 and proved the following  theorem  (for Bessel systems, not just orthogonal sets).
 \begin{theorem}\cite[Lemma 3.1.]{Lev18}\label{lev-1}
 Suppose that $\Lambda$ is an orthogonal set for $\mu$, where $\mu$ is a measure satisfying \eqref{lev-0} for some $\alpha$ with $0 \leqslant \alpha \leqslant d$. Then
$$
 \sup _{x \in \mathbb{R}^d} \#(\Lambda \cap B(x, R)) \leqslant C R^\alpha \quad(R\geqslant 1)
$$
 for a certain constant $C$ which does not depend on $R$.
 \end{theorem}
 Now, we can prove Theorem \ref{curve1} by using above the result.
 \begin{proof}[Proof of Theorem \ref{curve1}] According to Theorem \ref{lev-1}, our goal is to prove that there exists a constant $C>0$ and $R_0\ge 1$ such that for all $R\geq R_0\geq1$ $$ \int_{|t|<R}|\widehat{\rho}(t)|^2 d t\geq CR^{d-1}.$$
Let us write $B_R:=B_R(0)$, and define
\begin{align*}
  \psi_R(t) &:=\frac{1}{|B_{R/2}|}\,\big(\mathbf{1}_{B_{R/2}} * \mathbf{1}_{B_{R/2}}\big)(t)\\
   &=\frac{1}{|B_{R/2}(0)|}\int_{R^d} \mathbf{1}_{B_{R/2}(0)}(t-s)\,\mathbf{1}_{B_{R/2}(0)}(s)\,ds\\&=\frac{|B_{R/2}\cap (t+B_{R/2})|}{|B_{R/2}|}.
\end{align*}
For every $R>0$, we have the following three basic properties for $\psi_R$. 
\begin{enumerate}
\item $0\leq\psi_R\leq1$ and $\supp\psi_R \subset B_R$;
\item For all $x\in \mathbb{R}^d$ we have
\ifdefined\SMART
\begin{align*}
\hat{\psi}_R(x) &= \frac{1}{|B_{R/2}|}\,\Big| \hat{\mathbf{1}}_{B_{R/2}}(x)\Big|^2\\
 &= \frac{ 1 }{|B_{R/2}|}\cdot\left(\frac{R}{2 }\right)^{2d}\cdot\,\Big| \hat{\mathbf{1}}_{B_{1}}\left(\frac{R}{2}x\right)\Big|^2\\
&\ge 0;
\end{align*}
\else
$\hat{\psi}_R(x) = \frac{1}{|B_{R/2}|}\,\Big| \hat{\mathbf{1}}_{B_{R/2}}(x)\Big|^2
 = \frac{ 1 }{|B_{R/2}|}\cdot\left(\frac{R}{2 }\right)^{2d}\cdot\,\Big| \hat{\mathbf{1}}_{B_{1}}\left(\frac{R}{2}x\right)\Big|^2
\ge 0;$
\fi
\item Since $\hat{\mathbf{1}}_{B_{1}}$ is continuous and  $\hat{\mathbf{1}}_{B_{1}}(0)=|B_{1}|$, there exist constants $\varepsilon_0\in(0,1)$ and $\omega_d>0$ (depending only on $d$) such that  $\hat{\psi}_R(x) \ \ge \omega_d R^d$ for all $|x|\le \frac{\varepsilon_0}{R}$.
\end{enumerate}
Hence, we have
\begin{align*}
 \int_{\RR^d} \psi_R(t)\,|\hat{\rho}(t)|^2\,dt\;
  &= \iint \bigg(\int_{\RR^d} \psi_R(t)\, e^{-2\pi i\, t\cdot(x-y)}\,dt\bigg)\, d\rho(x)\,d\rho(y)\\&=\iint\hat{\psi}_R(x-y)\,d\rho(x)\,d\rho(y)
       \\&\ge\; \iint_{|x-y|\le\frac{\varepsilon_0}{R}} \hat{\psi}_R(x-y)\,d\rho(x)\,d\rho(y)
       \\&\ge \omega_d R^d\iint_{|x-y|\le \frac{\varepsilon_0}{R}} d\rho(x)\,d\rho(y).
\end{align*}
Moreover, since $\supp\psi_R \subset B_R$,
$$\int_{R^d} \psi_R(t)\,|\hat{\rho}(t)|^2\,dt=\;\int_{|t|<R} \psi_R(t)\,|\hat{\rho}(t)|^2\,dt\leq\;\int_{|t|<R} |\hat{\rho}(t)|^2\,dt.$$
For any curve $\Gamma_0$ in $\Gamma$ the restriction $\rho|_{\Gamma_0}$ is the arc-length measure on $\Gamma_0$.
 Write $\gamma_0:[0,L_0]\to\mathbb{R}^d$ for the arc-length parametrization of $\Gamma_0$ with $0<L_0<\infty$.  

We have
\begin{align*}\int_{|t|<R} |\hat{\rho}(t)|^2\,dt &\ge \omega_d R^d\iint_{|x-y|\le \frac{\varepsilon_0}{R}} d\rho(x)\,d\rho(y)\\
&\geq  
\omega_d R^d\,\iint \mathbf{1}_{\{|x-y|\le \frac{\varepsilon_0}{R}\}}\,\mathbf{1}_{\{x\in\gamma_{k_0},\,y\in\gamma_{k_0}\}}\,d\rho(x)\,d\rho(y)
\\&= \omega_d R^d\,  \int_{\gamma_{k_0}}\int_{\gamma_{k_0}} \mathbf{1}_{\{|x-y|\le \frac{\varepsilon_0}{R}\}}\, dx\, dy\\&=\omega_d R^d\, \int_0^{L_{k_0}}\int_0^{L_{k_0}} \mathbf{1}_{\{|\gamma_{k_0}(s)-\gamma_{k_0}(t)|\le\frac{\varepsilon_0}{R}\}}\,ds\,dt.
\end{align*}
Since $\gamma_{k_0} $ is arc-length parametrizations of curve segments, 
  we have  $|\gamma'_{k_0}(s)|\equiv1$ and $ |\gamma_{k_0}(s)-\gamma_{k_0}(t)|\le \frac{\varepsilon_0}{R}$ when $|s-t|\le \frac{\varepsilon_0}{R}$. By simple calculation, we can get
\ifdefined\SMART
\begin{align*}\label{eq:5.1}
\int_0^{L_{k_0}}\int_0^{L_{k_0}} \mathbf{1}_{\{|\gamma_{k_0}(s)-\gamma_{k_0}(t)|\le\frac{\varepsilon_0}{R}\}}\,ds\,dt &\ge\; \iint_{\substack{s,t\in[0,L_{k_0}]\\ |s-t|\le\frac{\varepsilon_0}{R}}} ds\,dt\\
&= \begin{cases}
L_{k_0}^2, & \frac{\varepsilon_0}{R}\geq L_{k_0},\\[6pt]
\frac{2\varepsilon_0L_{k_0}}{R}-(\frac{\varepsilon_0}{R})^2, & \frac{\varepsilon_0}{R}< L_{k_0}.
\end{cases}
\end{align*}
\else
\begin{equation*}\label{eq:5.1}
\int_0^{L_{k_0}}\int_0^{L_{k_0}} \mathbf{1}_{\{|\gamma_{k_0}(s)-\gamma_{k_0}(t)|\le\frac{\varepsilon_0}{R}\}}\,ds\,dt \ge\; \iint_{\substack{s,t\in[0,L_{k_0}]\\ |s-t|\le\frac{\varepsilon_0}{R}}} ds\,dt= \begin{cases}
L_{k_0}^2, & \frac{\varepsilon_0}{R}\geq L_{k_0},\\[6pt]
\frac{2\varepsilon_0L_{k_0}}{R}-(\frac{\varepsilon_0}{R})^2, & \frac{\varepsilon_0}{R}< L_{k_0}.
\end{cases} \end{equation*}
\fi
When $R>\frac{\varepsilon_0}{L_{k_0}}$, we can obtain  $$
\int_{|t|< R}|\widehat \rho(t)|^2\,dt
   \geq \omega_d \varepsilon_0 L_{k_0} R^{d-1}.$$
Taking $R_0=\max\{1, \frac{\varepsilon_0}{L_{k_0}}\}$ and $C=\omega_d \varepsilon_0 L_{k_0}$, we have $  \int_{|t|<R}|\widehat{\rho}(t)|^2 d t\geq CR^{d-1} $   for all $R\geq R_0 $. This completes the proof.
  
\end{proof}

\subsection{Ahlfors--David regular measures}\label{ss:shi}
 We present another approach to Theorem \ref{curve1} based on the second author's arguments in \cite{Shi21}. In fact, we prove more general results on Ahlfors--David regular measures in any dimension (not just for one-dimensional measures).

 We first recall several definitions and results in \cite{Shi21}.  We define the $n$-th dyadic partition of $\mathbb{R}$ by
$$
\mathcal{D}_n^{(1)}:=\left\{ \left[ \frac{k}{2^n}, \frac{k+1}{2^n} \right): k\in \mathbb{Z}  \right\},\ \ n\in\ZZ.
$$
The $n$-th dyadic partition of $\mathbb{R}^d$ is then defined by
$$
\mathcal{D}_n^{(d)}:=\left\{ I_1\times I_2\times \dots \times I_d: I_j\in \mathcal{D}_n^{(1)}  \right\}.
$$
If there is no confusion, we usually omit the superscript and write $\mathcal{D}_n$ for the $n$-th dyadic partition of $\mathbb{R}^d$.
Let $\mu\in \mathcal{P}(\mathbb{R}^d)$, where $\mathcal{P}(\mathbb{R}^d)$ stands for the space of probability measures on $\RR^d$. For a Borel set $K\subset \mathbb{R}^d$, we denote by 
$$
\mu_{K}(\cdot):=\mu(\cdot\cap K),
$$
the measure $\mu$ restricted on $K$. Moreover, if $K$ is a dyadic cube in $[0,1]^d$ with $\mu(K)>0$, we denote by
$$
\mu_{K}^{\Box}(\cdot):=\frac{1}{\mu(K)}(S_K)_*\mu_{K}(\cdot),
$$
where $S_K$ is the affine bijective map from $K$ to $[0,1]^d$ and $(S_K)_*\mu_{K}(\cdot)$ is the pushforward of $\mu_K$, i.e. the measure $\mu_{K}(S_K^{-1}(\cdot))$. Obviously, we have $\mu_{K}^{\Box}\in \mathcal{P}([0,1]^d)$.

\begin{lemma}[\cite{Shi21}, Lemma 4.2]\label{lem:fourier transform}
	Let $\mu\in \mathcal{P}(\mathbb{R}^d)$. Then for any $0<\epsilon<1$, there exists $\delta=\delta(\epsilon, \mu)>0$ such that for any $|\xi|<\delta$ and for any $D\in \mathcal{D}_n$ with $\mu(D)>0$, we have
	$$
	\left| \widehat{\mu_{D}^{\Box}}(\xi)  \right|>\epsilon.
	$$
\end{lemma}
For $\mu\in \mathcal{P}([0,1]^d)$ and $t,\lambda\in \RR^d$, we write the inner products
$$
\left\langle t, \lambda \right\rangle_{\mu}=\int_{[0,1]^d} e^{2\pi i (t-\lambda) \cdot x} d\mu(x).
$$

 We will only use the upper (Bessel-type) estimate from \cite[Lemma 4.1]{Shi21}. While \cite[Lemma 4.1]{Shi21} is proved in the frame spectral setting and yields two-sided bounds, our proof only requires the upper estimate (the Bessel side). We therefore state the result below in the form needed here, assuming only that 
$\Lambda$ is an orthogonal set for $\mu$.
 
\begin{lemma}[\cite{Shi21}, Lemma 4.1]\label{lem:change measure}
Let  $\mu\in \mathcal{P} (\RR^d) $ and $\Lambda\subset\RR^d$ be an orthogonal set for $\mu$, i.e.
    \begin{equation}\label{eq:change measure 1}
	  \sum_{\lambda\in \Lambda} \left|\left\langle t, \lambda \right\rangle_{\mu}  \right|^2\le  \mu(D).
	\end{equation}
    Let $n>0$. Then for any $D\in \mathcal{D}_n$ with $\mu(D)>0$, and for any $t\in \RR^d$, we have
	\begin{equation*}
    \sum_{\lambda\in \Lambda} \left|\left\langle \frac{1}{2^n}t, \frac{1}{2^n}\lambda \right\rangle_{\mu_{D}^{\Box}}  \right|^2\le \frac{1}{\mu(D)}.
	\end{equation*}
\end{lemma}

The following theorem follows the same proof as \cite[Theorem 1.2]{Shi21}. It also works for frame spectral measures, but we state it here for any measures in $\mathcal{P}(\mathbb{R}^d)$.
 
\begin{theorem}\label{thm:estimation}
   For any   measure $\mu  \in \mathcal{P}(\mathbb{R}^d)$  and   any orthogonal set $\Lambda$ for $\mu$, there exist constants $C$ and $\varrho$  such that  $t\in \mathbb{R}^d$ and $h>0$, we have
    \begin{equation}
	\# (\Lambda\cap B(t, h))
	 \le  C 2^{\sum_{D\in \mathcal{D}_n} -\mu(D)\log \mu(D)},
	 \end{equation}
     where $2^{n-1+\varrho}<h\le 2^{n+\varrho}$.
\end{theorem}
\begin{proof}
	 Notice that
	 \begin{equation}\label{eq:main thm 1}
	 2^{\sum_{D\in \mathcal{D}_n} -\mu(D)\log \mu(D)}=\prod_{D\in \mathcal{D}_n, \mu(D)>0} \mu(D)^{-\mu(D)}.
	 \end{equation}
	 Let $h>0$. Then there exists a positive integer $n_h$ such that $2^{n_h-1}< h\le 2^{n_h}$. Let $\epsilon>0$ and let $\delta=\delta(\epsilon, \mu)$, which is defined in Lemma \ref{lem:fourier transform}. Let $\varrho$ be the minimal integer such that $2^{-\varrho}<\delta$. For any $t\in \mathbb{R}^d$ and any $D\in \mathcal{D}_{n_h+\varrho}$, noticing that 
     $$
     \left|\frac{1}{2^{n_h+\varrho}}t-\frac{1}{2^{n_h+\varrho}}\lambda\right|<\delta, \text{ for any } \lambda\in \Lambda\cap B(t, h),
     $$
     by Lemma \ref{lem:fourier transform}, we have

\begin{equation}\label{eq:main thm 2}
	 \begin{split}
	 \epsilon^2 \cdot \#(\Lambda\cap B(t, h))
	 &\le \sum_{\lambda\in \Lambda\cap B(t, h)} \left|\left\langle \frac{1}{2^{n_h+\varrho}}t, \frac{1}{2^{n_h+\varrho}}\lambda \right\rangle_{\mu_{D}^{\Box}}  \right|^2\\
	 &\le \sum_{\lambda\in \Lambda} \left|\left\langle \frac{1}{2^{n_h+\varrho}}t, \frac{1}{2^{n_h+\varrho}}\lambda \right\rangle_{\mu_{D}^{\Box}}  \right|^2\\
 &\le \frac{1}{\mu({D})}.
	 \end{split}
	 \end{equation}
     The last inequality is due to   Lemma \ref{lem:change measure}. Since $\sum_{D\in \mathcal{D}_{n_h+\varrho}}\mu(D)=1$ and \eqref{eq:main thm 2} holds for all $D\in \mathcal{D}_{n_h+\varrho}$ with $\mu(D)>0$, we have
	 \begin{equation}\label{eq:main thm 3}
	 \begin{split}
\# (\Lambda\cap B(t, h))
	 &=\prod_{D\in \mathcal{D}_{n_h+\varrho}, \mu(D)>0}(\# (\Lambda\cap B(t, h)))^{\mu(D)} \\
	 &\le  \epsilon^{-2} \prod_{D\in \mathcal{D}_{n_h+\varrho}, \mu(D)>0} \mu(D)^{-\mu(D)}\\
  &=\epsilon^{-2} 2^{\sum_{D\in \mathcal{D}_{n_h+\varrho}} -\mu(D)\log \mu(D)}.
	 \end{split}
	 \end{equation}
	 The proof is complete.
     
\end{proof}

The measure $\mu$ is called \textbf{Ahlfors--David regular} (or simply \textbf{Ahlfors regular}) 
of dimension $s \geq 0$ if there exist constants $c, C > 0$ such that for every ball 
$B(x,r)$ with center $x \in \operatorname{supp}(\mu)$ and radius 
$0 < r < \operatorname{diam}(\operatorname{supp}(\mu))$, we have

\[
c\, r^s \leq \mu\bigl(B(x,r)\bigr) \leq C\, r^s .
\]

By applying Theorem \ref{thm:estimation}, we have the following theorem.
 
\begin{theorem}
For any orthogonal set $\Lambda$ of Ahlfors--David regular measures of dimension $s$, we have the estimate
$$\# (\Lambda\cap B(t, h)) = O(h^s), \text{ for } t\in \RR^d \text{ and } h>0. $$ 
In particular, we have the following estimates 
    \begin{enumerate}
        \item $\#(\Lambda\cap B(t, h)) = O(h)$ for any orthogonal set $\Lambda$ of any finite union of $C^1$ curves (for example, line segments).
        \item $\# (\Lambda\cap B(t, h)) = O(h^n)$ for any orthogonal set $\Lambda$ of any finite union of $n$-dimensional polyhedrons.
    \end{enumerate}
\end{theorem}
\begin{proof}
    We prove it for an Ahlfors--David regular measure $\mu$ of dimension $s\ge 0$. For any such measure there exist constants $c, C > 0$ such that for every $n$-th dyadic set
$D \subset \RR^d$ with center in $\supp \mu$, we have
\[
c\, 2^{-ns} \leq \mu\bigl(D\bigr) \leq C\, 2^{-ns} .
\]
By $\sum_{D\in \mathcal{D}_{n}}\mu(D)=1$, we have
$$
2^{\sum_{D\in \mathcal{D}_n} -\mu(D)\log \mu(D)} \le 2^{\sum_{D\in \mathcal{D}_n} -\mu(D)\log c2^{-ns}}=c2^{ns}\le \hat{C} h^s,
$$
for some constant $\hat{C}$. Then by Theorem \ref{thm:estimation}, we have 
$\# (\Lambda\cap B(t, h))=O(h^s)$.

\end{proof}

\bibliographystyle{alpha}

\end{document}